\documentclass{amsart}
\usepackage{amssymb} 
\usepackage{color}
\usepackage[all]{xy}
\usepackage{enumerate}

 \newtheorem{theorem}{Theorem}[section]

 \newtheorem{proposition}[theorem]{Proposition}
 \newtheorem{lemma}[theorem]{Lemma}
 \newtheorem{corollary}[theorem]{Corollary}
 \newtheorem{remark}[theorem]{Remark}
 \newtheorem{definition}[theorem]{Definition}
 \newtheorem{notation}[theorem]{Notation}
 \newtheorem{example}[theorem]{Example}
 \newtheorem{question}[theorem]{Question}

 \newtheorem{claim}[theorem]{Claim}
 \newtheorem{condition}[theorem]{Condition}
 
 \numberwithin{equation}{section}

 \def\subrel#1#2{\mathrel{\mathop{#2}\limits_{#1}}}

\def\cat{{\rm Cat}_{\infty}}
\def\catkappa{\cat^{(\kappa)}}

\def\overcat#1{\cat\mbox{$\scriptstyle /#1$}}

\def\wcat{\widehat{\rm Cat}_{\infty}}
\def\fin{{\rm Fin}_{\ast}}
\def\operad{{\rm Op}_{\infty}}
\def\operadgen{{\rm Op}_{\infty}^{\rm gen}}
\def\woperadgen{\widehat{\rm Op}_{\infty}^{\rm gen}}

\def\op#1{{\rm Op}_{\infty/#1}}
\def\wop#1{\widehat{\rm Op}_{\infty/#1}}

\def\alg{{\rm Alg}}

\def\map#1{{\rm Map}_{#1}}

\def\modenc#1{\mathrm{Mod}^{\mathbb{E}_n}_{#1}(\mathcal{C})}
\def\modemnc#1{\mathrm{Mod}^{\mathbb{E}_n/\mathbb{E}_{m+n}}_{#1}(\mathcal{C})}
\def\overlinemodemnc{\overline{\mathrm{Mod}}\,{}^{\opden/\opdemn}_{}(\mathcal{C})}
\def\widetildemodemnc{\widetilde{\mathrm{Mod}}\,{}^{\opden/\opdemn}_{}(\mathcal{C})}
\def\opden{\mathbb{E}_n}
\def\opdem{\mathbb{E}_m}
\def\opdemn{\mathbb{E}_{m+n}}
\def\algc#1{\mathrm{Alg}_{#1}(\mathcal{C})}
\def\palgc#1{{}^{\mathrm{p}}\mathrm{Alg}_{#1}(\mathcal{C})}
\def\poverlinealgc#1{{}^{\mathrm{p}}\overline{\mathrm{Alg}}_{#1}(\mathcal{C})}
\def\ko{\mathcal{K}_{\mathcal{O}}}
\def\kot{\ko^{\rm tr}}

\def\xko{{}_X\ko}
\def\xkot{\xko^{\rm tr}}
\def\xaoo{{}_X\mathcal{M}_{\mathcal{O}}^{\otimes}}

\def\xaot{\xao^{\rm tr}}
\def\xaoto{\xao^{{\rm tr},\otimes}}
\def\xabo{{}_X\widetilde{\mathcal{M}}_{\mathcal{O}}}
\def\xaboo{{}_X\widetilde{\mathcal{M}}_{\mathcal{O}}^{\otimes}}
\def\xao{{}_X\mathcal{M}_{\mathcal{O}}}

\begin{document}

\title
{Duoidal $\infty$-categories
of operadic modules}
\author{Takeshi Torii}
\address{Department of Mathematics, 
Okayama University,
Okayama 700--8530, Japan}
\email{torii@math.okayama-u.ac.jp}


\date{April 23, 2022 (version~1.0)}

\begin{abstract}
In this paper we study duoidal structures
on $\infty$-categories of operadic modules.
Let $\mathcal{O}^{\otimes}$ be a small coherent $\infty$-operad 
and let $\mathcal{P}^{\otimes}$ be an $\infty$-operad. 
If a $\mathcal{P}\otimes\mathcal{O}$-monoidal
$\infty$-category $\mathcal{C}^{\otimes}$
has a sufficient supply of colimits,
then we show that 
the $\infty$-category
${\rm Mod}_A^{\mathcal{O}}(\mathcal{C})$
of $\mathcal{O}$-$A$-modules in $\mathcal{C}^{\otimes}$
has a structure of $(\mathcal{P},\mathcal{O})$-duoidal 
$\infty$-category
for any $\mathcal{P}\otimes\mathcal{O}$-algebra
object $A$.

\end{abstract}

\maketitle

\section{Introduction}

A duoidal category is a category equipped with
two monoidal structures in which 
one is lax monoidal with respect to the other.
In \cite{Torii1,Torii3} we have introduced 
generalizations of duoidal categories 
in the setting of $\infty$-categories.
The goal of this paper is
to show that $\infty$-categories of operadic modules
have duoidal structures.

The notion of duoidal category was introduced by
Aguiar-Mahajan~\cite{Aguiar-Mahajan} 
by the name of $2$-monoidal category.
There is a $2$-category 
${\rm Mon}^{\rm oplax}({\rm Cat})$ of monoidal categories,
oplax monoidal functors, and natural transformations
between them.
Note that ${\rm Mon}^{\rm oplax}({\rm Cat})$ 
is a monoidal $2$-category under Cartesian product.
A duoidal category is identified with 
a pseudomonoid in the monoidal $2$-category 
${\rm Mon}^{\rm oplax}({\rm Cat})$.

We can consider the $\infty$-category
${\rm Mon}_{\mathcal{O}}^{\rm oplax}(\cat)$
of $\mathcal{O}$-monoidal $\infty$-categories
and oplax $\mathcal{O}$-monoidal functors
for an $\infty$-operad $\mathcal{O}^{\otimes}$.
Since it has finite products,
${\rm Mon}_{\mathcal{O}}^{\rm oplax}(\cat)$
is a Cartesian symmetric monoidal $\infty$-category.
For an $\infty$-operad $\mathcal{P}^{\otimes}$,
we say that a $\mathcal{P}$-monoid object
in the $\infty$-category
${\rm Mon}_{\mathcal{O}}^{\rm oplax}(\cat)$
is a $(\mathcal{P},\mathcal{O})$-duoidal
$\infty$-category.

Now we assume that $\mathcal{O}^{\otimes}$ is coherent.
For an $\mathcal{O}$-monoidal
$\infty$-category $\mathcal{C}^{\otimes}$ 
which has a sufficient supply of colimits, 
Lurie~\cite{Lurie2} constructed an $\mathcal{O}$-monoidal structure
on the $\infty$-category of $\mathcal{O}$-$A$-modules
in $\mathcal{C}^{\otimes}$
for each $\mathcal{O}$-algebra object $A$.
The main theorem in this paper
is to extend this $\mathcal{O}$-monoidal structure
to a $(\mathcal{P},\mathcal{O})$-duoidal structure.

\begin{theorem}[{cf.~Theorem~\ref{thm:main-functor}}]
Let $\kappa$ be an uncountable regular cardinal
and let $\mathcal{O}^{\otimes}$ be an essentially
$\kappa$-small coherent $\infty$-operad.
Let $\mathcal{P}^{\otimes}$ be an $\infty$-operad 
and let $\mathcal{C}^{\otimes}$ be 
a $\mathcal{P}\otimes\mathcal{O}$-monoidal $\infty$-category 
which is compatible with $\kappa$-small colimits.
Then the $\infty$-category
${\rm Mod}_A^{\mathcal{O}}(\mathcal{C})$
of $\mathcal{O}$-$A$-modules in $\mathcal{C}^{\otimes}$
has a structure of a $(\mathcal{P},\mathcal{O})$-duoidal
$\infty$-category
for any $\mathcal{P}\otimes\mathcal{O}$-algebra object $A$.
\end{theorem}

The important case is
when $(\mathcal{P}^{\otimes},\mathcal{O}^{\otimes})
=(\mathbb{E}_m^{\otimes},\mathbb{E}_n^{\otimes})$,
where $\mathbb{E}_k$ is the little $k$-cubes operad,
and $\mathcal{C}^{\otimes}$ is a presentable 
symmetric monoidal $\infty$-category.
In this case we have the following corollary.

\begin{corollary}[{cf.~Theorem~\ref{thm:Em-En-case}}]
Let $(\mathcal{P}^{\otimes},\mathcal{O}^{\otimes})=
(\mathbb{E}_m^{\otimes},\mathbb{E}_n^{\otimes})$ and
let $\mathcal{C}^{\otimes}$ be a presentable
symmetric monoidal $\infty$-category.
Then
the $\infty$-category ${\rm Mod}_A^{\mathbb{E}_n}(\mathcal{C})$
of $\mathbb{E}_n$-$A$-modules in $\mathcal{C}^{\otimes}$
has a structure of 
an $(\mathbb{E}_m,\mathbb{E}_n)$-duoidal $\infty$-category
for any $\mathbb{E}_{m+n}$-algebra object $A$.
\end{corollary}

The organization of this paper is as follows:
In \S\ref{section:operadic-modules}
we study $\infty$-categories of operadic modules.
First, we recall the construction
of a map of generalized $\infty$-operads 
which encodes $\infty$-categories of operadic modules
and restriction functors.
Then we give a description
of left adjoints to the restriction functors.
In \S\ref{section:adjointable-diagrams}
we study adjointable diagrams 
of (op)lax $\mathcal{O}$-monoidal functors
between $\mathcal{O}$-monoidal $\infty$-categories.
In \S\ref{section:mixed-fibration-operadic-module}
we construct a map of generalized $\infty$-operads
which encodes $\infty$-categories of 
operadic modules and left adjoint oplax monoidal functors
by using the results in \S\ref{section:adjointable-diagrams}.
Then we study a functoriality of this construction.
In \S\ref{section:duoidal-operadic-module}
we will construct duoidal structures
on $\infty$-categories of operadic modules 
and prove the main theorem (Theorem~\ref{thm:main-functor}). 
In \S\ref{section:em-en-duoidal}
we consider the case in which 
$(\mathcal{P}^{\otimes},\mathcal{O}^{\otimes})=
(\mathbb{E}_m^{\otimes},\mathbb{E}_n^{\otimes})$
and $\mathcal{C}^{\otimes}$ is a presentable
symmetric monoidal $\infty$-category.

\begin{notation}\rm
For $\infty$-operads
$\mathcal{O}^{\otimes}$ and $\mathcal{P}^{\otimes}$,
we denote by $\mathcal{O}^{\otimes}\otimes\mathcal{P}^{\otimes}
=(\mathcal{O}\otimes \mathcal{P})^{\otimes}$
the Boardman-Vogt tensor product of $\infty$-operads
(cf.~\cite[\S2.2.5]{Lurie2}).
For an $\infty$-operad $\mathcal{O}^{\otimes}$,
we denote by
${\rm Mon}_{\mathcal{O}}^{\rm lax}(\cat)$
the $\infty$-category
of small $\mathcal{O}$-monoidal
$\infty$-categories and lax $\mathcal{O}$-monoidal
functors.
We also denote by
${\rm Mon}_{\mathcal{O}}^{\rm oplax}(\cat)$
the $\infty$-category
of small $\mathcal{O}$-monoidal
$\infty$-categories and oplax $\mathcal{O}$-monoidal
functors.
\end{notation}

\section{Operadic modules}
\label{section:operadic-modules}

In this section we study $\infty$-categories of operadic modules.
Let $p: \mathcal{O}^{\otimes}\to {\rm Fin}_*$ 
be a coherent $\infty$-operad
and 
let $q: \mathcal{C}^{\otimes}\to \mathcal{O}^{\otimes}$ 
be a map of $\infty$-operads.
In \S\ref{subsection:infinity-cat-operadic-modules}
we recall a generalized $\infty$-operad
${\rm Mod}^{\mathcal{O}}(\mathcal{C})^{\otimes}$
constructed by Lurie,
in which the underlying $\infty$-category
${\rm Mod}^{\mathcal{O}}(\mathcal{C})$ consists
of pairs $(A,M)$ of an $\mathcal{O}$-algebra object $A$
and an $\mathcal{O}$-$A$-module $M$
in $\mathcal{C}^{\otimes}$.
We would like to construct a 
free functor
${\rm Alg}_{/\mathcal{O}}(\mathcal{C})\times
\mathcal{C}^{\otimes}_X\to {\rm Mod}^{\mathcal{O}}(\mathcal{C})^{\otimes}_X$
for $X\in\mathcal{O}$,
which is left adjoint to the forgetful functor.
For this purpose,
in \S\ref{subsection:operad-xaoo} 
and \S\ref{subsection:operad-xaoto}
we introduce $\infty$-operads $\xaoo$ 
and $\xaoto$
such that
the $\infty$-categories of $\xao$-algebras 
and of $\xaot$-algebras in $\mathcal{C}^{\otimes}$
are equivalent to 
${\rm Mod}^{\mathcal{O}}(\mathcal{C})^{\otimes}_X$
and 
${\rm Alg}_{/\mathcal{O}}(\mathcal{C})\times
\mathcal{C}^{\otimes}_X$,
respectively.
In \S\ref{subsection:free-operadic-modules}
we construct a left adjoint $f_{!X}$ to 
the restriction functor $f^*_X$
for a map $f: A\to B$ in ${\rm Alg}_{/\mathcal{O}}(\mathcal{C})$.
We also study the monad $\mathbf{T}_f$
associated to the adjunction
$(f_{!X},f^*_X)$ and give a description
of $\mathbf{T}_f(M)$ as a colimit
of certain diagram for 
$M\in {\rm Mod}_A^{\mathcal{O}}(\mathcal{C})^{\otimes}_X$.

\subsection{$\infty$-categories of operadic modules}
\label{subsection:infinity-cat-operadic-modules}

Let $p: \mathcal{O}^{\otimes}\to {\rm Fin}_*$ 
be a coherent $\infty$-operad
and 
let $q: \mathcal{C}^{\otimes}\to \mathcal{O}^{\otimes}$ 
be a map of $\infty$-operads.
In \cite[\S3.3.3]{Lurie2}
Lurie introduced a generalized $\infty$-operad
${\rm Mod}^{\mathcal{O}}(\mathcal{C})^{\otimes}$,
where the underlying $\infty$-category
${\rm Mod}^{\mathcal{O}}(\mathcal{C})$
consists of pairs $(A,M)$,
where $A$ is an $\mathcal{O}$-algebra and
$M$ is an $\mathcal{O}$-$A$-module in $\mathcal{C}^{\otimes}$.
In this subsection
we recall the construction
of ${\rm Mod}^{\mathcal{O}}(\mathcal{C})^{\otimes}$
and a map 
$(\Psi,\sigma): {\rm Mod}^{\mathcal{O}}(\mathcal{C})^{\otimes}
\to {\rm Alg}_{/\mathcal{O}}(\mathcal{C})\times
\mathcal{O}^{\otimes}$
of generalized $\infty$-operads.

First,
we recall the construction 
of the generalized $\infty$-operad
${\rm Mod}^{\mathcal{O}}(\mathcal{C})^{\otimes}$.
Let $\ko$
be a full subcategory of ${\rm Fun}([1],\mathcal{O}^{\otimes})$
spanned by semi-inert morphisms in $\mathcal{O}^{\otimes}$
(see \cite[Definition~3.3.1.1]{Lurie2} for the definition of
semi-inert morphisms).
We have the projections
${\rm ev}_0,{\rm ev}_1:\ko\to \mathcal{O}^{\otimes}$
given by evaluation at $0,1\in [1]$,
respectively.
A morphism in $\ko$ is said to be inert
if the images under ${\rm ev}_0$ and ${\rm ev}_1$
are inert morphisms in $\mathcal{O}^{\otimes}$. 

We have an $\infty$-category
$\widetilde{\rm Mod}{}^{\mathcal{O}}(\mathcal{C})^{\otimes}$
equipped with a map
$\widetilde{\rm Mod}{}^{\mathcal{O}}(\mathcal{C})^{\otimes}
\to \mathcal{O}^{\otimes}$
so that there is an equivalence
\[ {\rm Map}_{{\rm Cat}_{\infty/\mathcal{O}^{\otimes}}}
    (\mathcal{X},\widetilde{\rm Mod}{}^{\mathcal{O}}(\mathcal{C})^{\otimes})
   \simeq
    {\rm Map}_{{\rm Cat}_{\infty/\mathcal{O}^{\otimes}}}
   (\mathcal{X}\times_{\mathcal{O}^{\otimes},{\rm ev}_0}
                \ko,\mathcal{C}^{\otimes})\]
for any map $\mathcal{X}\to\mathcal{O}^{\otimes}$
of $\infty$-categories.
An object of $\widetilde{\rm Mod}{}^{\mathcal{O}}
(\mathcal{C})^{\otimes}$
over $Y\in\mathcal{O}^{\otimes}$
corresponds to a functor
$\{Y\}\times_{\mathcal{O}^{\otimes},{\rm ev}_0}\ko\to
\mathcal{C}^{\otimes}$ over $\mathcal{O}^{\otimes}$.
We let
\[ {\rm Mod}^{\mathcal{O}}(\mathcal{C})^{\otimes} \]
be the full subcategory of
$\widetilde{\rm Mod}{}^{\mathcal{O}}(\mathcal{C})^{\otimes}$
spanned by those functors 
$\{Y\}\times_{\mathcal{O}^{\rm op},{\rm ev}_0}\ko\to
\mathcal{C}^{\otimes}$ 
which preserve inert morphisms.
The induced map ${\rm Mod}^{\mathcal{O}}(\mathcal{C})^{\otimes}\to
\mathcal{O}^{\otimes}$
is a map of generalized $\infty$-operads
by \cite[Proposition~3.3.3.10]{Lurie2}.

Next, we recall a map
${\rm Mod}^{\mathcal{O}}(\mathcal{C})^{\otimes}
   \to
   {\rm Alg}_{/\mathcal{O}}(\mathcal{C})
   \times\mathcal{O}^{\otimes}$
of generalized $\infty$-operads
and its properties.
Let $\ko^0$ be the full subcategory
of $\ko$ spanned by null morphisms
(see \cite[Definition~3.3.1.1]{Lurie2} 
for the definition of null morphisms).
By \cite[Remark~3.3.3.16]{Lurie2},
the inclusion
$\ko^0\hookrightarrow\ko$
induces a map
\[ (\Phi,\sigma):   
   {\rm Mod}^{\mathcal{O}}(\mathcal{C})^{\otimes}
   \longrightarrow 
   {\rm Alg}_{/\mathcal{O}}(\mathcal{C})
   \times\mathcal{O}^{\otimes}\]
of generalized $\infty$-operads.
By taking fibers at $A\in {\rm Alg}_{/\mathcal{O}}(\mathcal{C})$,
we obtain a functor
\[ \sigma_A: {\rm Mod}_A^{\mathcal{O}}(\mathcal{C})^{\otimes}
             \longrightarrow
             \mathcal{O}^{\otimes},\]
which is a map of $\infty$-operads
by \cite[Theorem~3.3.3.9]{Lurie2}.
By \cite[Corollary~3.4.3.4]{Lurie2},
the map $\Phi: {\rm Mod}^{\mathcal{O}}(\mathcal{C})^{\otimes}
\to {\rm Alg}_{/\mathcal{O}}(\mathcal{C})$
is a Cartesian fibration,
and the induced functor 
\[ f^*: {\rm Mod}_B^{\mathcal{O}}(\mathcal{C})^{\otimes}
          \longrightarrow
          {\rm Mod}_A^{\mathcal{O}}(\mathcal{C})^{\otimes}, \]
is a map of $\infty$-operads 
over $\mathcal{O}^{\otimes}$
for any map
$f: A\to B$ in ${\rm Alg}_{/\mathcal{O}}(\mathcal{C})$.

If $\mathcal{C}^{\otimes}$ has a sufficient supply
of colimits, then
$\sigma_A: {\rm Mod}_A^{\mathcal{O}}(\mathcal{C})^{\otimes}\to
\mathcal{O}^{\otimes}$
is a coCartesian fibration of $\infty$-operads
for any $A\in {\rm Alg}_{/\mathcal{O}}(\mathcal{C})$.
Let $\kappa$ be an uncountable regular cardinal.
We assume that $\mathcal{O}^{\otimes}$
is an essentially $\kappa$-small coherent $\infty$-operad
and that $\mathcal{C}^{\otimes}$
is an $\mathcal{O}$-monoidal $\infty$-category
which is compatible with $\kappa$-small colimits
(see \cite[Definition~3.1.1.18 and Variant~3.1.1.19 ]{Lurie2}
for the definition of $\mathcal{O}$-monoidal $\infty$-categories
compatible with $\kappa$-small colimits). 
Then $\sigma_A: {\rm Mod}_A^{\mathcal{O}}(\mathcal{C})^{\otimes}
\to\mathcal{O}^{\otimes}$
is an $\mathcal{O}$-monoidal $\infty$-category
which is compatible with $\kappa$-small colimits
by \cite[Corollary~3.4.4.6]{Lurie2}.


\subsection{The $\infty$-operad $\xaoo$}
\label{subsection:operad-xaoo}

We would like to construct a 
free functor
${\rm Alg}_{/\mathcal{O}}(\mathcal{C})\times
\mathcal{C}^{\otimes}_X\to {\rm Mod}^{\mathcal{O}}(\mathcal{C})^{\otimes}_X$
for $X\in\mathcal{O}$,
which is left adjoint to the forgetful functor.
For this purpose,
in this subsection
we introduce an $\infty$-operad $\xaoo$ such that
the $\infty$-category ${\rm Alg}_{\xao/\mathcal{O}}(\mathcal{C})$
of $\xao$-algebras in $\mathcal{C}^{\otimes}$
is equivalent to 
${\rm Mod}^{\mathcal{O}}(\mathcal{C})^{\otimes}_X$.

For $X\in\mathcal{O}$,
we set 
\[ \xko=
   \{X\}\times_{\mathcal{O}^{\otimes},{\rm ev}_0}\ko. \]
The $\infty$-category $\xko$
is equipped with a map 
$\xko\to\mathcal{O}^{\otimes}$
induced by ${\rm ev}_1$.
We notice that 
the $\infty$-category ${\rm Mod}^{\mathcal{O}}(\mathcal{C})_X^{\otimes}$
is a full subcategory of 
${\rm Fun}_{\mathcal{O}^{\otimes}}(\xko,\mathcal{C}^{\otimes})$
spanned by those functors
which preserve inert morphisms.

%

Let ${\rm Triv}^{\otimes}$ be the trivial $\infty$-operad,
which is the subcategory of ${\rm Fin}_*$ spanned
by inert morphisms
(\cite[Example~2.1.1.20]{Lurie2}).
The inclusion map
$\{\langle 1\rangle\}\hookrightarrow
{\rm Triv}^{\otimes}$
induces an equivalence
${\rm Alg}_{\rm Triv}(\mathcal{O})\simeq\mathcal{O}$
by \cite[Remark~2.1.3.6]{Lurie2}.
Hence we have a map of $\infty$-operads
$X: {\rm Triv}^{\otimes}\to \mathcal{O}^{\otimes}$
for $X\in\mathcal{O}$.
We set
\[ \xaboo=
   {\rm Triv}^{\otimes}\times_{X,\mathcal{O}^{\otimes},{\rm ev}_0}
   \ko.\]
We define 
\[ \xaoo \]
to be the full subcategory of 
$\xaboo$
spanned by those vertices
which correspond to semi-inert morphisms 
$\oplus_m X\to Y_1\oplus\cdots\oplus Y_n$
in $\mathcal{O}^{\otimes}$
such that the image 
$\langle m\rangle \to \langle n\rangle$
in ${\rm Fin}_*$ is an order-preserving injection.

We will show that 
the composite map
$r: \xaoo\stackrel{{\rm ev}_1}{\to} \mathcal{O}^{\otimes}
   \stackrel{p}{\to} {\rm Fin}_*$
exhibits $\xaoo$ as an $\infty$-operad,
and that 
${\rm ev}_1: \xaoo\to\mathcal{O}^{\otimes}$
is a map of $\infty$-operads.

\begin{proposition}\label{prop:xaoo-infinity-operad}
The map 
$r: \xaoo\to {\rm Fin}_*$
exhibits $\xaoo$ as an $\infty$-operad.
\end{proposition}

\proof
Let $\alpha: \oplus_mX\to Y\simeq Y_1\oplus\cdots\oplus Y_n$
be an object of $\xaoo$
over $\underline{\alpha}:\langle m\rangle \to 
\langle n\rangle$ in ${\rm Fin}_*$.
We show that there is an $r$-coCartesian 
morphism $\rho^i_!: \alpha\to \beta$ over 
$\rho^i:\langle n\rangle \to\langle 1\rangle$. 
We decompose the composite
$\overline{\rho}^i_!\circ \alpha: \oplus_m X\to Y\to Y_i$
as an inert morphism $\oplus_m X\to Z$
followed by an active morphism $\beta: Z\to Y_i$,
where $\overline{\rho}^i_!: Y\to Y_i$
is a $p$-coCartesian morphism in $\mathcal{O}^{\otimes}$
over $\rho^i$.
This determines a morphism 
$\rho^i_!:\alpha\to \beta$ in $\xaoo$
over $\rho^i$.
We can verify that $\rho^i_!: \alpha\to \beta$ 
is an $r$-coCartesian morphism.

By using the fact that
$p:\mathcal{O}^{\otimes}\to{\rm Fin}_*$
is an $\infty$-operad,
we can verify that
the $r$-coCartesian morphisms $\rho^i_!$ for $1\le i\le n$
induces an equivalence 
$\prod_{1\le i\le n}\rho^i_!:
(\xaoo)_{\langle n\rangle}\to 
\prod_{1\le i\le n}(\xaoo)_{\langle 1\rangle}$
of $\infty$-categories
for any $n\ge 0$, and 
that the family 
$\{\rho^i_!: \alpha\to \rho^i_!(\alpha)\}_{1\le i\le n}$
of morphisms is an $r$-limit diagram
in $\xaoo$.
\qed

\begin{corollary}
The map ${\rm ev}_1: \xaoo\to \mathcal{O}^{\otimes}$
is a map of $\infty$-operads.
\end{corollary}

\proof
Let $\alpha: \oplus_m X\to Y\simeq Y_1\oplus\cdots\oplus Y_n$
be an object of $\xaoo$,
and let $\rho^i: \langle n\rangle \to\langle 1\rangle$
be the inert morphism in ${\rm Fin}_*$.
By the proof of Proposition~\ref{prop:xaoo-infinity-operad},
we have ${\rm ev}_1(\rho^i_!: \alpha\to \rho^i_!(\alpha))
=(\overline{\rho}^i_!:
Y\to Y_i)$.
Hence we see that ${\rm ev}_1$ preserves inert morphisms.
\qed

\bigskip

Now, we will introduce an $\infty$-category
$\xao^{\le 1}$
which is equivalent to $\xko$
and show that $\xao^{\le 1}$
is an approximation to 
the $\infty$-operad $\xaoo$.

We let $\xabo^{\le 1}$
be the full subcategory
of $\xaboo$
spanned by those vertices
which correspond to maps
$\oplus_mX\to Y_1\oplus\cdots\oplus Y_n$
for $0\le m\le 1$.
Note that 
$\xko$ is a full subcategory of $\xabo^{\le 1}$.
We define
\[ \xao^{\le 1} \]
to be the full subcategory
of $\xaoo$ spanned by those vertices
which correspond to maps
$\oplus_mX\to Y_1\oplus\cdots\oplus Y_n$
for $0\le m\le 1$.
Note that $\xao^{\le 1}$ is a full subcategory
of $\xabo^{\le 1}$.

We consider a right Kan extension
of the identity functor $\xao^{\le 1}\to\xao^{\le 1}$ 
along the inclusion functor
$\xao^{\le 1}\hookrightarrow \xabo^{\le 1}$:
\[ \xymatrix{
   \xao^{\le 1}\ar[r]\ar[d] & \xao^{\le 1}\\
   \xabo^{\le 1}\ar@{..>}[ur]_R&\\
}\]

Let $\alpha: X\to Y\simeq Y_1\oplus\cdots\oplus Y_n$ be an 
object of $\xabo^{\le 1}$,
which is not contained in $\xao^{\le 1}$.
Note that $\alpha$ is a null morphism
(see \cite[Definition~3.3.1.1]{Lurie2} 
for the definition of null morphisms).
We observe that
the $\infty$-category 
$\xao^{\le 1}\times_{\xabo^{\le 1}}\xabo^{\le 1}{}_{/\alpha}$
has an initial object $\langle 0\rangle\to Y$.
Hence there exists a right
Kan extension $R: \xabo^{\le 1}\to\xao^{\le 1}$ 
by \cite[Lemma~4.3.2.13]{Lurie1}.

By restricting $R$ to $\xko$,
we obtain a functor
\[ R: \xko\to \xao^{\le 1}. \]
We will show that the functor $R$ gives an equivalence
of $\infty$-categories.

\begin{lemma}\label{lemma:R-equivalence}
The functor 
$R: \xko\to \xao^{\le 1}$
is an equivalence of $\infty$-categories.
\end{lemma}

\proof
By the same argument as above,
we see that there exists a left Kan extension
$L: \xabo^{\le 1}\to\xko$
of the identity functor $\xko\to\xko$
along the inclusion functor 
$\xko\to \xabo^{\le 1}$.
By restricting $L$ to the full subcategory $\xao^{\le 1}$
of $\xabo^{\le 1}$,
we obtain a functor $L: \xao^{\le 1}\to \xko$. 
We can easily verify that the pair
$(L,R)$ of functors gives
an equivalence of $\infty$-categories.
\qed

\bigskip

Next, we will show that
${\rm Mod}^{\mathcal{O}}(\mathcal{C})^{\otimes}_X$
is equivalent to
the $\infty$-category ${\rm Alg}_{\xao/\mathcal{O}}(\mathcal{C})$
of $\xao$-algebra objects in $\mathcal{C}^{\otimes}$.

\begin{lemma}
The inclusion functor
$i: \xao^{\le 1}\hookrightarrow \xaoo$
is an approximation to 
the $\infty$-operad $\xaoo$.
\end{lemma}

\proof
We will verify the conditions
in \cite[Definition~2.3.3.6]{Lurie2}.
Since $\xao^{\le 1}$ is a full subcategory
of $\xaoo$,
we see that condition~(1) holds.
Let $\alpha: \oplus_mX\to Y$ be an object of $\xao^{\le 1}$
and let $\beta: \oplus_{n}X\to Z$ be an object of $\xaoo$.
We suppose that there exists an active morphism 
$\phi: \beta\to \alpha$.
Since $\xao^{\le 1}$ is a full subcategory
of $\xaoo$,
in order to show that condition~(2) holds,
it suffices to show that $\beta$ is an object
of $\xao^{\le 1}$.
The existence of the active morphism $\phi$
implies that $m=n$.
Thus, we see that $\beta\in\xao^{\le 1}$.
\qed

\bigskip

By \cite[Definition~2.3.3.20]{Lurie2},
we have an $\infty$-category
${\rm Alg}_{\xao^{\le 1}}(\mathcal{X})$
for an $\infty$-operad $\mathcal{X}^{\otimes}\to{\rm Fin}_*$,
which is a full subcategory of
${\rm Fun}_{{\rm Fin}_*}(\xao^{\le 1},\mathcal{X}^{\otimes})$
spanned by those functors which preserve
inert morphisms.
The map $q:\mathcal{C}^{\otimes}\to\mathcal{O}^{\otimes}$
of $\infty$-operads induces a map
$q_*: {\rm Alg}_{\xao^{\le 1}}(\mathcal{C})\to
{\rm Alg}_{\xao^{\le 1}}(\mathcal{O})$
of $\infty$-categories. 
We define
\[ {\rm Alg}_{\xao^{\le 1}/\mathcal{O}}(\mathcal{C}) \]
to be the fiber of the map $q_*$ 
at ${\rm ev}_1\circ i: \xao^{\le 1}\to \xaoo\to 
\mathcal{O}^{\otimes}$.

\begin{lemma}
\label{lemma:xao-le-1-algebra}
The inclusion functor
$i: \xao^{\le 1}\hookrightarrow \xaoo$
induces an equivalence
\[ i^*: {\rm Alg}_{\xao/\mathcal{O}}(\mathcal{C})
   \stackrel{\simeq}{\longrightarrow}
   {\rm Alg}_{\xao^{\le 1}/\mathcal{O}}(\mathcal{C}).\]
\end{lemma}

\proof
We have a commutative diagram
\[ \begin{array}{ccc}
    {\rm Alg}_{\xao}(\mathcal{C})
    & \stackrel{q_*}{\longrightarrow} &
    {\rm Alg}_{\xao}(\mathcal{O})\\
    \mbox{$\scriptstyle i^*$}
    \bigg\downarrow 
    \phantom{\mbox{$\scriptstyle i^*$}}
    & & 
    \phantom{\mbox{$\scriptstyle i^*$}}
    \bigg\downarrow
    \mbox{$\scriptstyle i^*$} \\
    {\rm Alg}_{\xao^{\le 1}}(\mathcal{C})
    & \stackrel{q_*}{\longrightarrow} &
    {\rm Alg}_{\xao^{\le 1}}(\mathcal{O})\\
   \end{array}\]
of $\infty$-categories.
The lemma follows from
the fact that 
the vertical arrows are equivalences
by \cite[Theorem~2.3.3.23]{Lurie2}.
\qed

\begin{proposition}\label{cor:equivalence-Alg-le-1}
The composite
$i\circ R: \xko\stackrel{}{\to}
\xao^{\le 1}\stackrel{}{\to}
\xaoo$ 
induces an equivalence
\[ {\rm Alg}_{\xao/\mathcal{O}}(\mathcal{C})
   \stackrel{\simeq}{\longrightarrow}
   {\rm Mod}^{\mathcal{O}}(\mathcal{C})^{\otimes}_X. \]
\end{proposition}

\proof
Since $R$ is an equivalence  
by Lemma~\ref{lemma:R-equivalence}
and preserves inert morphisms,
it induces an equivalence
$R^*: {\rm Mod}^{\mathcal{O}}(\mathcal{C})^{\otimes}_X
\stackrel{\simeq}{\to}{\rm Alg}_{\xao^{\le 1}/\mathcal{O}}(\mathcal{C})$.
The proposition follows from
Lemma~\ref{lemma:xao-le-1-algebra}.
\qed


\subsection{The operad $\xaoto$}
\label{subsection:operad-xaoto}

We would like to construct a 
free functor
${\rm Alg}_{/\mathcal{O}}(\mathcal{C})\times
\mathcal{C}^{\otimes}_X\to {\rm Mod}^{\mathcal{O}}(\mathcal{C})^{\otimes}_X$
for $X\in\mathcal{O}$,
which is left adjoint to the forgetful functor.
For this purpose,
in this subsection
we introduce an $\infty$-operad 
$\xaoto$
such that
the $\infty$-category ${\rm Alg}_{\xaoto/\mathcal{O}}(\mathcal{C})$
of $\xaot$-algebras in $\mathcal{C}^{\otimes}$
is equivalent to 
${\rm Alg}_{/\mathcal{O}}(\mathcal{C})\times
\mathcal{C}^{\otimes}_X$.

First,
we introduce an $\infty$-operad $\xaoto$
such that the $\infty$-category 
${\rm Alg}_{\xaot/\mathcal{O}}(\mathcal{C})$
of $\xaot$-algebras in $\mathcal{C}^{\otimes}$
is equivalent to
${\rm Mod}^{\mathcal{O}}(\mathcal{C})_X^{\rm tr}$.
We define 
\[ \xaoto \]
to be a subcategory of $\xaoo$
as follows:
The objects of $\xaoto$
are the same as those of $\xaoo$.
A morphism  
\[ \begin{array}{ccc}
   \oplus_m X
   &\stackrel{\alpha}{\longrightarrow}&
   Y_1\oplus\cdots\oplus Y_n\\
   \bigg\downarrow & & 
   \phantom{\mbox{$\scriptstyle \gamma$}}
   \bigg\downarrow
   \mbox{$\scriptstyle \gamma$} \\
   \oplus_{m'}X
   &\stackrel{\beta}{\longrightarrow}&
   Z_1\oplus\cdots\oplus Z_{n'}\\
  \end{array} \]
in $\xaoo$ is a morphism 
in $\xaoto$ if and only if 
the image
\[ \begin{array}{ccc}
   \langle m\rangle & 
   \stackrel{\underline{\alpha}}{\longrightarrow} &
   \langle n\rangle\\
   \bigg\downarrow & &
   \phantom{\mbox{$\scriptstyle \underline{\gamma}$}}
   \bigg\downarrow
   \mbox{$\scriptstyle \underline{\gamma}$}\\
   \langle m'\rangle &
   \stackrel{\underline{\beta}}{\longrightarrow} &
   \langle n'\rangle\\ 
  \end{array}\]
in ${\rm Fin}_*$ satisfies 
the condition that the cardinality of the set
$\underline{\gamma}^{-1}(\underline{\beta}(i))$ is just
one for each $i\in \langle m'\rangle^{\circ}$.
We can verify that 
$\xaoto$ supports
a structure of $\infty$-operad
such that the inclusion functor
$\xaoto\hookrightarrow \xaoo$ is a map
of $\infty$-operads.

Next, we will show that
the $\infty$-category
${\rm Alg}_{\xaot/\mathcal{O}}(\mathcal{C})$
of $\xaot$-algebras in $\mathcal{C}^{\otimes}$
is equivalent to
${\rm Alg}_{/\mathcal{O}}(\mathcal{C})\times
\mathcal{C}^{\otimes}_X$.
For this purpose,
we construct
an equivalence
$\mathcal{O}^{\otimes}\boxplus {\rm Triv}^{\otimes}
\stackrel{\simeq}{\to} \xao^{0,\otimes}$
of $\infty$-operads,
where the left hand side is
a coproduct of $\infty$-operads. 

We define 
\[ \xao^{0,\otimes}=
   \{\langle 0\rangle\}\times_{{\rm Triv}^{\otimes}}\xaoo.
 \]
We can verify that 
the map $r: \xaoo\to {\rm Fin}_*$ induces a map
$r: \xao^{0,\otimes}\to {\rm Fin}_*$
which exhibits $\xao^{0,\otimes}$ as 
an $\infty$-operad. 
Furthermore,
the restriction of the map
${\rm ev}_1: \xaoo\to\mathcal{O}^{\otimes}$ induces
an equivalence
${\rm ev}_1: \xao^{0,\otimes}\stackrel{\simeq}{\longrightarrow}
   \mathcal{O}^{\otimes}$
of $\infty$-operads.

We notice that $\xao^{0,\otimes}$
is a full subcategory of $\xaoto$
and that the inclusion functor
$\xao^{0,\otimes}\hookrightarrow \xaoto$
is a map of $\infty$-operads. 
Thus, we obtain a map
$\mathcal{O}^{\otimes}\to \xaoto$
of $\infty$-operads
since $\xao^{0,\otimes}$
is equivalent to $\mathcal{O}^{\otimes}$.
Furthermore, since there is an equivalence
${\rm Alg}_{\rm Triv}(\xao^{\rm tr})\simeq
(\xaoto)_{\langle 1\rangle}$
and $({\rm id}_X: X\to X)$
is an object of $(\xaoto)_{\langle 1\rangle}$,
there is a map
${\rm Triv}^{\otimes}\to \xaoto$
of $\infty$-operads
corresponding to ${\rm id}_X$.
Combining these two maps,
we obtain a map
\[ \mathcal{O}^{\otimes}\boxplus 
   {\rm Triv}^{\otimes}\longrightarrow
   \xaoto \]
of $\infty$-operads.

By the definition of coproducts of $\infty$-operads
(\cite[\S2.2.3]{Lurie2}),
we easily obtain the following lemma.

\begin{lemma}
The map 
$\mathcal{O}^{\otimes}\boxplus 
 {\rm Triv}^{\otimes}\stackrel{\simeq}{\to}\xaoto$
is an equivalence of $\infty$-operads.
\end{lemma}

\begin{corollary}\label{cor:xaot-algebra-description}
There is an equivalence
\[ {\rm Alg}_{\xaot/\mathcal{O}}(\mathcal{C})
   \stackrel{\simeq}{\longrightarrow}
   {\rm Alg}_{/\mathcal{O}}(\mathcal{C})\times\mathcal{C}^{\otimes}_X.
 \]
\end{corollary}

\begin{remark}\rm
The inclusion functors
$\xao^{0,\otimes}\hookrightarrow
 \xaoto\hookrightarrow \xaoo$
induce the following commutative diagram
\[ \begin{array}{ccccc}
    {\rm Alg}_{\xao/\mathcal{O}}(\mathcal{C})
    &\longrightarrow&
    {\rm Alg}_{\xaot/\mathcal{O}}(\mathcal{C})
    & \longrightarrow &
    {\rm Alg}_{\xao^{0}/\mathcal{O}}(\mathcal{C})\\
    \phantom{\mbox{$\scriptstyle \simeq$}}
    \bigg\downarrow
    \mbox{$\scriptstyle \simeq$}
    & & 
    \phantom{\mbox{$\scriptstyle \simeq$}}  
    \bigg\downarrow
    \mbox{$\scriptstyle \simeq$}
    & & 
    \phantom{\mbox{$\scriptstyle \simeq$}}
    \bigg\downarrow
    \mbox{$\scriptstyle \simeq$} \\[3mm]
    {\rm Mod}^{\mathcal{O}}(\mathcal{C})^{\otimes}_X
    &\longrightarrow&
    {\rm Alg}_{/\mathcal{O}}(\mathcal{C})\times\mathcal{C}^{\otimes}_X
    &\longrightarrow&
    {\rm Alg}_{/\mathcal{O}}(\mathcal{C}).\\
   \end{array}\]
\end{remark}

\subsection{Free operadic modules}
\label{subsection:free-operadic-modules}

For a map $f: A\to B$ in ${\rm Alg}_{/\mathcal{O}}(\mathcal{C})$
and $X\in \mathcal{O}$,
we have the restriction functor 
$f^*_X: {\rm Mod}^{\mathcal{O}}_B(\mathcal{C})^{\otimes}_X
\to {\rm Mod}^{\mathcal{O}}_A(\mathcal{C})^{\otimes}_X$.
In this subsection
we first construct a left adjoint $f_{!X}$ to $f^*_X$
under some conditions. 
After that,
we study the monad $\mathbf{T}_f$
associated to the adjunction
$(f_{!X},f^*_X)$ 
and describe ${\mathbf T}_f(M)$
as a colimit of certain diagram
for $M\in {\rm Mod}_A^{\mathcal{O}}(\mathcal{C})^{\otimes}_X$.

Let $\kappa$ be an uncountable regular cardinal.
In this subsection
we assume that $\mathcal{O}^{\otimes}$ is an
essentially $\kappa$-small 
coherent $\infty$-operad,
and that $q: \mathcal{C}^{\otimes}\to \mathcal{O}^{\otimes}$
is an $\mathcal{O}$-monoidal $\infty$-category
which is compatible with $\kappa$-small colimits
in the sense of 
\cite[Definition~3.1.1.18 and Variant~3.1.1.19]{Lurie2}.

We set 
$\mathcal{D}^{\otimes}={\rm Mod}_A^{\mathcal{O}}
(\mathcal{C})^{\otimes}$
for simplicity.
By \cite[Corollary~3.4.4.6]{Lurie2},
the map 
$\mathcal{D}^{\otimes}\to \mathcal{O}^{\otimes}$
is an $\mathcal{O}$-monoidal $\infty$-category
compatible with $\kappa$-small colimits.
We can regard $A$ as an object of 
$\mathcal{D}^{\otimes}$,
and it is a unit object of the $\mathcal{O}$-monoidal structure.
Furthermore,
there is an equivalence
${\rm Alg}_{/\mathcal{O}}(\mathcal{D})
\simeq {\rm Alg}_{/\mathcal{O}}(\mathcal{C})_{A/}$
of $\infty$-categories by \cite[Corollary~3.4.1.7]{Lurie2}.
Thus, we can regard $f$ as an $\mathcal{O}$-algebra
object of $\mathcal{D}^{\otimes}$.
By \cite[Corollary~3.4.1.9]{Lurie2},
there is an equivalence
${\rm Mod}_{f}^{\mathcal{O}}
(\mathcal{D})^{\otimes}
\stackrel{\simeq}{\to}
{\rm Mod}_B^{\mathcal{O}}(\mathcal{C})^{\otimes}$
of $\mathcal{O}$-monoidal $\infty$-categories.

Let $j: \xaoto\to \xaoo$ be the inclusion map,
which is a map of $\infty$-operads.
The map $j$ induces a functor
\[ j^*: {\rm Mod}^{\mathcal{O}}(\mathcal{D})^{\otimes}_X
        \longrightarrow
        {\rm Alg}_{/\mathcal{O}}(\mathcal{D})
        \times\mathcal{D}^{\otimes}_X.\]
By \cite[Corollary3.1.3.5]{Lurie2},
the functor $j^*$ admits a left adjoint 
\[ j_!:  {\rm Alg}_{/\mathcal{O}}(\mathcal{D})\times 
       \mathcal{D}^{\otimes}_X
        \longrightarrow
        {\rm Mod}^{\mathcal{O}}(\mathcal{D})^{\otimes}_X, \]
which is obtained by
the operadic left Kan extension along $j$.
Note that $j_!$ makes the following diagram commute
\begin{align}\label{align:j-adjunction-triangle}
    \xymatrix{
    {\rm Alg}_{/\mathcal{O}}(\mathcal{D})\times
    \mathcal{D}^{\otimes}_X
    \ar[rr]^{j_!}\ar[dr] &&
    {\rm Mod}^{\mathcal{O}}(\mathcal{D})^{\otimes}_X
    \ar[dl] \\
    &{\rm Alg}_{/\mathcal{O}}(\mathcal{D}).&\\
    }
\end{align}

By taking fibers at 
$f\in {\rm Alg}_{/\mathcal{O}}(\mathcal{D})$
in diagram~(\ref{align:j-adjunction-triangle}),
we obtain a functor
$f_{!X}: \mathcal{D}^{\otimes}_X
           \to
           {\rm Mod}_{B}^{\mathcal{O}}
(\mathcal{D})^{\otimes}_X$,
which is identified with a functor
\[ f_{!X}: {\rm Mod}_A^{\mathcal{O}}(\mathcal{C})^{\otimes}_X
           \longrightarrow
           {\rm Mod}_B^{\mathcal{O}}(\mathcal{C})^{\otimes}_X. \]

\begin{proposition}\label{prop:f-adjunction-general}
The functor $f_{!X}$ is a left adjoint 
to $f_X^*$.
\end{proposition}

\proof
For any 
$(f,M)\in {\rm Alg}_{/\mathcal{O}}(\mathcal{D})\times
\mathcal{D}^{\otimes}_X$ and 
$(f,N)\in{\rm Mod}^{\mathcal{O}}(\mathcal{D})^{\otimes}_X$,
we have a map
\[    {\rm Map}_{{\rm Mod}^{\mathcal{O}}(\mathcal{D})^{\otimes}_X}
   (j_!(f,M),(f,N))
   \longrightarrow
   {\rm Map}_{{\rm Alg}_{/\mathcal{O}}(\mathcal{D})\times\mathcal{D}^{\otimes}_X}
   ((f,M),j^*(f,N)) \]
of mapping spaces
over ${\rm Map}_{{\rm Alg}_{/\mathcal{O}}(\mathcal{D})}(f,f)$,
which is an equivalence 
since $(j_!,j^*)$ is an adjoint pair.
By taking fibers at ${\rm id}_f\in
{\rm Map}_{{\rm Alg}_{/\mathcal{O}}(\mathcal{D})}(f,f)$,
we obtain a natural equivalence
\[ {\rm Map}_{{\rm Mod}_f^{\mathcal{O}}(\mathcal{D})^{\otimes}_X}
   (f_{!X}(M),N)
   \stackrel{\simeq}{\longrightarrow}
   {\rm Map}_{\mathcal{D}^{\otimes}_X}(M,f_X^*(N))\]
of mapping spaces,
which completes the proof.
\qed

\bigskip



We let 
\[ \mathbf{T}_f=f^*_X\circ f_{!X}:
    {\rm Mod}_A^{\mathcal{O}}(\mathcal{C})^{\otimes}_X
    \longrightarrow
    {\rm Mod}_A^{\mathcal{O}}(\mathcal{C})^{\otimes}_X\]
be the monad associated 
to the adjunction $(f_{!X},f^*_X)$.
Next, we study the monad $\mathbf{T}_f$
and describe ${\mathbf T}_f(M)$
as a colimit of certain diagram.

We set
\[ I(\mathcal{O})=(\xaoto)_{\rm act}\times_{(\xaoo)_{\rm act}}
   {(\xaoo)_{\rm act}}_{/{\rm id}_X}.\]
For $(f,M)\in {\rm Alg}_{/\mathcal{O}}(\mathcal{D})
\times \mathcal{D}^{\otimes}_X$,
we will construct a functor
\[ \overline{D}(f,M): I(\mathcal{O})^{\triangleright}
   \longrightarrow \mathcal{D}^{\otimes}_X
   \simeq {\rm Mod}_A^{\mathcal{O}}(\mathcal{C})^{\otimes}_X,\]
which is a colimit diagram and 
carries the cone point to 
$\mathbf{T}_f(M)$.

By Corollary~\ref{cor:xaot-algebra-description},
a pair $(f,M)\in 
{\rm Alg}_{/\mathcal{O}}(\mathcal{D})\times
\mathcal{D}^{\otimes}_X$
determines a map
\[ (f,M): \xaoto\longrightarrow \mathcal{D}^{\otimes} \]
over $\mathcal{O}^{\otimes}$.
By the operadic left Kan extension $j_!$,
we have an operadic $q$-colimit diagram 
\[ \xymatrix{
    I(\mathcal{O})\ar[r]\ar[d]& 
    (\xaoto)_{\rm act}\ar[r]&
    \xaoto\ar[r]^{(f,M)}&
    \mathcal{D}^{\otimes}\ar[d]^q\\
    I(\mathcal{O})^{\triangleright}\ar[r]\ar[urrr]&
    \left({(\xaoo)_{\rm act}}_{/{\rm id}_X}\right)^{\triangleright}\ar[r]&
    \xaoo\ar[r]_{{\rm ev}_1}& \mathcal{O}^{\otimes}\\ 
}\]
by \cite[Proposition~3.1.3.3]{Lurie2},
where the bottom arrows
carries the cone point of $I(\mathcal{O})^{\triangleright}$
to $X\in\mathcal{O}$.
By the coCartesian pushforward
of the map $I(\mathcal{O})^{\triangleright}\to \mathcal{D}^{\otimes}$,
we obtain a functor
$\overline{D}(f,M): I(\mathcal{O})^{\triangleright}
   \to \mathcal{D}^{\otimes}_X$,
which is a colimit diagram
and carries the cone point to 
$\mathbf{T}_f(M)$.

By the functoriality of construction,
we obtain the following proposition.

\begin{proposition}\label{prop:functor:alg-C-fun-I-C}
We have a functor
\[ D: {\rm Alg}_{/\mathcal{O}}(\mathcal{C})_{A/}\times
      {\rm Mod}_A^{\mathcal{O}}(\mathcal{C})^{\otimes}_X
      \longrightarrow 
      {\rm Fun}(I(\mathcal{O}),
      {\rm Mod}_A^{\mathcal{O}}(\mathcal{C})^{\otimes}_X) \]
which assigns to a pair 
$(f,M)\in {\rm Alg}_{/\mathcal{O}}(\mathcal{C})_{A/}
\times {\rm Mod}_A(\mathcal{C})^{\otimes}_X$
a diagram 
$D(f,M): I(\mathcal{O})
           \to {\rm Mod}_A^{\mathcal{O}}
            (\mathcal{C})^{\otimes}_X$
such that 
\[ \mathbf{T}_f(M)\simeq\ \subrel{I(\mathcal{O})}{\rm colim}
            D(f,M). \]
\end{proposition}

\section{Adjointable diagrams
of monoidal $\infty$-categories}
\label{section:adjointable-diagrams}

In this section 
we introduce $\infty$-categories
$\mathsf{Mod}_{\mathcal{O}}^{\rm oplax, L}(\cat)^{\rm RAd}$
and 
$\mathsf{Mod}_{\mathcal{O}}^{\rm lax, R}(\cat)^{\rm LAd}$,
and study their properties.
In \S\ref{subsection:equiv-diagram-infty-categories}
we construct 
$\mathsf{Mod}_{\mathcal{O}}^{\rm oplax, L}(\cat)^{\rm RAd}$
equipped with a Cartesian fibration
$\mathsf{Mod}_{\mathcal{O}}^{\rm oplax, L}(\cat)^{\rm RAd}\to\cat$
encoding diagrams of 
$\mathcal{O}$-monoidal $\infty$-categories,
left adjoint oplax $\mathcal{O}$-monoidal functors,
and 
right adjointable commutative squares.
We also construct 
$\mathsf{Mod}_{\mathcal{O}}^{\rm lax, R}(\cat)^{\rm LAd}$
equipped with a Cartesian fibration
$\mathsf{Mod}_{\mathcal{O}}^{\rm lax, R}(\cat)^{\rm LAd}\to\cat$
encoding diagrams of 
$\mathcal{O}$-monoidal $\infty$-categories,
right adjoint lax $\mathcal{O}$-monoidal functors,
and 
left adjointable commutative squares.
We show that there is an equivalence
of $\infty$-categories
between $\mathsf{Mod}_{\mathcal{O}}^{\rm lax, R}(\cat)^{\rm LAd}$
and 
$\mathsf{Mod}_{\mathcal{O}}^{\rm oplax, L}(\cat)^{\rm RAd}$
by taking left adjoints to right adjoint 
lax $\mathcal{O}$-monoidal functors.
In \S\ref{subsection:mon-LAd-sub-Op}
we show that 
$\mathsf{Mon}_{\mathcal{O}}^{\rm lax,R}(\cat)^{\rm LAd}$ 
is equivalent to a subcategory of
${\rm Op}_{\infty/\mathcal{O}^{\otimes}}^{\rm gen}$,
where ${\rm Op}_{\infty}^{\rm gen}$
is the $\infty$-category of generalized $\infty$-operads.
Dually, in \S\ref{subsection:mon-RAd-opposite-op}
we show that 
$\mathsf{Mon}_{\mathcal{O}}^{\rm oplax,L}(\cat)^{\rm RAd}$
is equivalent to a subcategory of
${\rm Op}_{\infty/\mathcal{O}^{\otimes,{\rm op}}}^{\rm gen,\vee}$.

\subsection{An equivalence between
$\mathsf{Mod}_{\mathcal{O}}^{\rm lax, R}(\cat)^{\rm LAd}$
and
$\mathsf{Mod}_{\mathcal{O}}^{\rm oplax, L}(\cat)^{\rm RAd}$}
\label{subsection:equiv-diagram-infty-categories}

In this subsection
we construct Cartesian fibrations 
$\mathsf{Mod}_{\mathcal{O}}^{\rm oplax, L}(\cat)^{\rm RAd}\to\cat$
and
$\mathsf{Mod}_{\mathcal{O}}^{\rm lax, R}(\cat)^{\rm LAd}\to\cat$.
By taking left adjoints to right adjoint 
lax $\mathcal{O}$-monoidal functors,
we show that there is an equivalence
of $\infty$-categories
between $\mathsf{Mod}_{\mathcal{O}}^{\rm lax, R}(\cat)^{\rm LAd}$
and 
$\mathsf{Mod}_{\mathcal{O}}^{\rm oplax, L}(\cat)^{\rm RAd}$.

We have a wide subcategory
\[ {\rm Mon}_{\mathcal{O}}^{\rm lax, R}(\cat) \]
of ${\rm Mon}_{\mathcal{O}}^{\rm lax}(\cat)$
spanned by those lax $\mathcal{O}$-monoidal 
functors $f^*$ such that
$f^*_X$ is right adjoint for each $X\in\mathcal{O}$.
For an $\infty$-category $S$,
we denote by
\[ {\rm Fun}^{\rm LAd}(S,{\rm Mon}_{\mathcal{O}}^{\rm lax,R}
   (\cat)) \]
the $\infty$-category whose objects are
functors $F: S\to {\rm Mon}_{\mathcal{O}}^{\rm lax,R}(\cat)$
and whose functors are natural transformations
$\alpha: F\to G$ such that
$\alpha(s): F(s)\to G(s)$ is strong $\mathcal{O}$-monoidal functors
for each $s\in S$, and that
the following commutative diagram
\[ \begin{array}{ccc}
    F(s)_X &  \longrightarrow & F(s')_X\\
    \bigg\downarrow && \bigg\downarrow \\
    G(s)_X &  \longrightarrow & G(s')_X\\ 
   \end{array}\]
is left adjointable
for any morphism $s\to s'$ in $S$ 
and any $X\in\mathcal{O}$
(see \cite[Definition~4.7.4.13]{Lurie2}
for the definition of left adjointable diagrams).

We also have a wide subcategory
\[ {\rm Mon}_{\mathcal{O}}^{\rm oplax,L}(\cat) \]
of ${\rm Mon}_{\mathcal{O}}^{\rm oplax}(\cat)$
spanned by
those oplax $\mathcal{O}$-monoidal functors
$f_!$ such that $f_{!X}$ is left adjoint
for each $X\in\mathcal{O}$.
We denote by
\[ {\rm Fun}^{\rm RAd}(S,{\rm Mon}_{\mathcal{O}}^{\rm oplax,L}
   (\cat)) \]
the $\infty$-category whose objects are
functors $F: S\to {\rm Mon}_{\mathcal{O}}^{\rm oplax,L}(\cat)$
and whose functors are natural transformations
$\alpha: F\to G$ such that
$\alpha(s): F(s)\to G(s)$ is strong $\mathcal{O}$-monoidal functors
for each $s\in S$, and that
the following commutative diagram
\[ \begin{array}{ccc}
    F(s)_X &  \longrightarrow & F(s')_X\\
    \bigg\downarrow && \bigg\downarrow \\
    G(s)_X &  \longrightarrow & G(s')_X\\ 
   \end{array}\]
is right adjointable
for any morphism $s\to s'$ in $S$ 
and any $X\in\mathcal{O}$
(see \cite[Definition~4.7.4.13]{Lurie2}
for the definition of right adjointable diagrams).

First,
we shall show that there is a natural equivalence
between
${\rm Fun}^{\rm LAd}(S,{\rm Mon}_{\mathcal{O}}^{\rm lax,R}
    (\cat))$ and
${\rm Fun}^{\rm RAd}(S^{\rm op},{\rm Mon}_{\mathcal{O}}^{\rm oplax,L}
    (\cat))$
by taking left adjoints to right adjoint
lax $\mathcal{O}$-monoidal functors.

\begin{proposition}\label{prop:monoidal-adjointable-equivalence}
For any $\infty$-category $S$,
there is a natural equivalence
\[ {\rm Fun}^{\rm LAd}(S,{\rm Mon}_{\mathcal{O}}^{\rm lax,R}
    (\cat))\simeq
   {\rm Fun}^{\rm RAd}(S^{\rm op},{\rm Mon}_{\mathcal{O}}^{\rm oplax,L}
    (\cat)) \]
of $\infty$-categories.
\end{proposition}

\proof
It suffices to show that there is a natural equivalence
\[ {\rm Map}_{\cat}([n],
    {\rm Fun}^{\rm LAd}(S,{\rm Mon}_{\mathcal{O}}^{\rm lax,R}
     (\cat))) 
    \simeq
    {\rm Map}_{\cat}([n],
    {\rm Fun}^{\rm RAd}(S^{\rm op},{\rm Mon}_{\mathcal{O}}^{\rm oplax,L}
     (\cat))) \]
for any $[n]\in\Delta^{\rm op}$.

We let 
$({\rm Cat}^{\rm coc}_{\infty/\mathcal{B}})_{\rm lax}^{\rm radj}$
be a subcategory of ${\rm Cat}_{\infty/\mathcal{B}}$
whose objects are coCartesian fibrations and
whose morphisms are $\mathcal{B}$-parametrized right adjoints,
and let 
$({\rm Cat}^{\rm cart}_{\infty/\mathcal{B}})_{\rm oplax}^{\rm ladj}$
be a subcategory of ${\rm Cat}_{\infty/\mathcal{B}}$
whose objects are Cartesian fibrations and
whose morphisms are $\mathcal{B}$-parametrized left adjoints
(see \cite[Definition~2.1]{Haugseng2}
for the definition of 
$\mathcal{B}$-parametrized left and right adjoints).
There are inclusions
\[ \begin{array}{rcl}
   {\rm Map}_{\cat}([n],
    {\rm Fun}^{\rm LAd}(S,{\rm Mon}_{\mathcal{O}}^{\rm lax,R}
     (\cat)))&\hookrightarrow&
   {\rm Map}_{\cat}(S,
     (\cat^{\rm coc}{}_{/[n]\times\mathcal{O}^{\otimes}}
      )_{\rm lax}^{\rm radj}),\\[2mm]
   {\rm Map}_{\cat}([n],
    {\rm Fun}^{\rm RAd}(S^{\rm op},
    {\rm Mon}_{\mathcal{O}}^{\rm oplax,L}
     (\cat)))&\hookrightarrow&
    {\rm Map}_{\cat}(S^{\rm op},
     (\cat^{\rm cart}{}_{/[n]^{\rm op}\times\mathcal{O}^{\otimes,{\rm op}}}
      )_{\rm oplax}^{\rm ladj})\\
   \end{array}\]
of mapping spaces.
By \cite[Theorem~2.2]{Haugseng2},
we have an equivalence
\[  {\rm Map}_{\cat}(S,
     (\cat^{\rm coc}{}_{/[n]\times\mathcal{O}^{\otimes}}
      )_{\rm lax}^{\rm radj})
    \simeq 
    {\rm Map}_{\cat}(S^{\rm op},
     (\cat^{\rm cart}{}_{/[n]^{\rm op}\times\mathcal{O}^{\otimes,{\rm op}}}
      )_{\rm oplax}^{\rm ladj}), \]
which restricts to the desired equivalence.
\qed

\bigskip

We define
\[ \mathsf{Mon}_{\mathcal{O}}^{\rm lax, R}(\cat)^{\rm LAd}
   \longrightarrow \cat \]
to be a Cartesian fibration which is obtained
by unstraightening
of the functor $\cat^{\rm op}\to\wcat$
given by
$S\mapsto {\rm Fun}^{\rm LAd}
(S,{\rm Mon}_{\mathcal{O}}^{\rm lax,R}(\cat))$.
We also define
\[ \mathsf{Mon}_{\mathcal{O}}^{\rm oplax, L}(\cat)^{\rm RAd}
   \longrightarrow \cat \]
to be a Cartesian fibration which is obtained
by unstraightening
of the functor $\cat^{\rm op}\to\wcat$
given by
$S\mapsto {\rm Fun}^{\rm RAd}
(S,{\rm Mon}_{\mathcal{O}}^{\rm oplax,L}(\cat))$.
By Proposition~\ref{prop:monoidal-adjointable-equivalence},
we obtain the following theorem.

\begin{theorem}\label{thm:lax-oplax-monoidal-adjointable-eq}
There is an equivalence
\[ \mathsf{Mon}_{\mathcal{O}}^{\rm lax,R}(\cat)^{\rm LAd}
   \stackrel{\simeq}{\longrightarrow}
   \mathsf{Mon}_{\mathcal{O}}^{\rm oplax,L}(\cat)^{\rm RAd} \]
of $\infty$-categories,
which fits into the following commutative diagram
\[ 
    \xymatrix{
     \mathsf{Mon}_{\mathcal{O}}^{\rm lax,R}(\cat)^{\rm LAd}
     \ar[r]^-{\simeq}\ar[d]&
     \mathsf{Mon}_{\mathcal{O}}^{\rm oplax, L}(\cat)^{\rm RAd}
     \ar[d]\\
     \cat \ar[r]^{(-)^{\rm op}}&
     \cat.\\
    }\]
\end{theorem}

\begin{remark}\rm
We can construct a double $\infty$-category
\[ \mathbf{Mon}_{\mathcal{O}}^{\rm lax,R}
   (\cat)^{\rm LAd}\]
whose objects are $\mathcal{O}$-monoidal $\infty$-categories,
whose horizontal $1$-morphisms are
strong $\mathcal{O}$-monoidal functors,
whose vertical $1$-morphisms are
right adjoint lax $\mathcal{O}$-monoidal functors, and
whose $2$-morphisms are commutative squares
which are left adjointable for each $X\in\mathcal{O}$.

We can also construct a double $\infty$-category
\[ \mathbf{Mon}_{\mathcal{O}}^{\rm oplax,L}
   (\cat)^{\rm RAd}\]
whose objects are $\mathcal{O}$-monoidal $\infty$-categories,
whose horizontal $1$-morphisms are
strong $\mathcal{O}$-monoidal functors,
whose vertical $1$-morphisms are
left adjoint oplax $\mathcal{O}$-monoidal functors, and
whose $2$-morphisms are commutative squares
which are right adjointable for each $X\in\mathcal{O}$.

There is an equivalence
\[ \mathbf{Mon}_{\mathcal{O}}^{\rm lax,R}(\cat)^{\rm LAd}
   \simeq
   (\mathbf{Mon}_{\mathcal{O}}^{\rm oplax,L}(\cat)^{\rm RAd}
    )^{\mbox{\scriptsize v-op}},\]
where the right hand side is 
the double $\infty$-category 
obtained from $\mathbf{Mon}_{\mathcal{O}}^{\rm oplax,L}(\cat)^{\rm RAd}$
by reversing the vertical direction.
\end{remark}

\subsection{$\mathsf{Mon}_{\mathcal{O}}^{\rm lax,R}(\cat)^{\rm LAd}$
and ${\rm Op}_{\infty/\mathcal{O}^{\otimes}}^{\rm gen}$}
\label{subsection:mon-LAd-sub-Op}

In this subsection
we show that 
$\mathsf{Mon}_{\mathcal{O}}^{\rm lax, R}(\cat)^{\rm LAd}$
is equivalent to a subcategory of
${\rm Op}_{\infty{/\mathcal{O}^{\otimes}}}^{\rm gen}$,
where ${\rm Op}_{\infty}^{\rm gen}$
is the $\infty$-category of generalized $\infty$-operads.
For this purpose,
we study a functor 
\[ (-)_{\langle 0\rangle}: 
   {\rm Op}_{\infty/\mathcal{O}^{\otimes}}^{\rm gen}\longrightarrow \cat \]
which assigns 
the fiber $\mathcal{E}^{\otimes}_{\langle 0\rangle}$
at $\langle 0\rangle\in\mathcal{O}^{\otimes}$
to a map $\mathcal{E}^{\otimes}\to\mathcal{O}^{\otimes}$
of generalized $\infty$-operads.

\begin{lemma}\label{lemma:0-fiber-cartesian-1}
The functor $(-)_{\langle 0\rangle}: 
{\rm Op}_{\infty/\mathcal{O}^{\otimes}}^{\rm gen}\to\cat$
is a Cartesian fibration.
\end{lemma}

\proof
Let
$F: {\rm Op}_{\infty}^{\rm gen}\to \cat$
be the functor which associates
the fiber $\mathcal{X}^{\otimes}_{\langle 0\rangle}$
at $\langle 0\rangle \in{\rm Fin}_*$
to a generalized $\infty$-operad 
$\mathcal{X}^{\otimes}\to{\rm Fun}_*$.
By \cite[Proposition~2.3.2.9]{Lurie2},
the functor $F$
has a right adjoint
$G: \cat\to{\rm Op}_{\infty}^{\rm gen}$
which is given by
$C\mapsto C\times{\rm Fin}_*$.

Let $\mathcal{E}^{\otimes}\to \mathcal{O}^{\otimes}$
be a map of generalized $\infty$-operads.
We set $S=\mathcal{E}^{\otimes}_{\langle 0\rangle}$.
By the unit of the adjunction
$(F,G)$,
we have a map $\mathcal{E}^{\otimes}\to 
S\times{\rm Fin}_*$
of generalized $\infty$-operads.
This induces a map 
$\mathcal{E}^{\otimes}\to 
(S\times{\rm Fin}_*)\times_{{\rm Fin}_*} \mathcal{O}^{\otimes}
\simeq S\times\mathcal{O}^{\otimes}$ 
of generalized $\infty$-operads.

For any functor $f: T\to S$ in $\cat$,
we can verify that the projection
$T\times_S\mathcal{E}^{\otimes}\to \mathcal{E}^{\otimes}$
is a Cartesian morphism in ${\rm Op}_{\infty/\mathcal{O}^{\otimes}}^{\rm gen}$
covering $f$.
Hence the map $(-)_{\langle 0\rangle}:
{\rm Op}_{\infty/\mathcal{O}^{\otimes}}^{\rm gen}\to\cat$ 
is a Cartesian fibration.
\qed

\bigskip

We recall that
$\mathsf{Mon}_{\mathcal{O}}^{\rm lax, R}(\cat)^{\rm LAd}\to\cat$
is a Cartesian fibration 
which is associated to 
the functor $\cat^{\rm op}\to\wcat$
given by
$S\mapsto {\rm Fun}^{\rm LAd}(S,
{\rm Mon}_{\mathcal{O}}^{\rm lax,R}(\cat))$.
By the definition of
the $\infty$-category
${\rm Fun}^{\rm LAd}(S,
{\rm Mon}_{\mathcal{O}}^{\rm lax,R}(\cat))$,
there is a natural map
${\rm Fun}^{\rm LAd}(S,
{\rm Mon}_{\mathcal{O}}^{\rm lax,R}(\cat))\to
{\rm Fun}(S,{\rm Op}_{\infty/\mathcal{O}^{\otimes}})$.

By \cite[\S3]{Torii1},
we can identify 
the $\infty$-category
${\rm Fun}(S,{\rm Op}_{\infty}{}_{/\mathcal{O}^{\otimes}})$
with a subcategory 
of the slice category $\cat{}_{/S\times\mathcal{O}^{\otimes}}$
for any $\infty$-category $S$.
For an object of 
${\rm Fun}(S,{\rm Op}_{\infty}{}_{/\mathcal{O}^{\otimes}})$,
the corresponding object 
is a map
\[ f: \mathcal{E}^{\otimes}\to S\times \mathcal{O}^{\otimes} \]
such that 

\begin{itemize}

\item 
the map $f_S: \mathcal{E}^{\otimes}\to S$
is a coCartesian fibration,
and $f$ preserves coCartesian morphisms,
where $f_S$ is the composite of $f$
with the projection $S\times\mathcal{O}^{\otimes}\to S$,

\item
the restriction $f_s: \mathcal{E}^{\otimes}_s\to \mathcal{O}^{\otimes}$
is a map of $\infty$-operads for each $s\in S$, and

\item
the induced map $\mathcal{E}^{\otimes}_s\to
\mathcal{E}^{\otimes}_{s'}$ over $\mathcal{O}^{\otimes}$
preserves inert morphisms for each
morphism $s\to s'$ in $S$.

\end{itemize}

We can easily verify the following lemma.

\begin{lemma}\label{lemma:map-S-op-generalized-operad}
The map $f: \mathcal{E}^{\otimes}\to S\times \mathcal{O}^{\otimes}$
is a morphism of ${\rm Op}_{\infty/\mathcal{O}^{\otimes}}^{\rm gen}$
such that $\mathcal{E}^{\otimes}_{\langle 0\rangle}\to S$
is an equivalence of $\infty$-categories.
\end{lemma}

By Lemma~\ref{lemma:map-S-op-generalized-operad},
we can identify ${\rm Fun}(S,{\rm Op}_{\infty/\mathcal{O}^{\otimes}})$
with a subcategory of 
${\rm Op}_{\infty/\mathcal{O}^{\otimes}}^{\rm gen}\times_{\cat}
\{S\}$.
Hence we obtain 
a natural transformation
\[ {\rm Fun}^{\rm LAd}
  (-,{\rm Mon}_{\mathcal{O}}^{\rm lax,R}(\cat))\to
  {\rm Op}_{\infty/\mathcal{O}^{\otimes}}^{\rm gen}\times_{\cat}\{-\} \]
in which
${\rm Fun}^{\rm LAd}(S,{\rm Mon}_{\mathcal{O}}^{\rm lax,R}(\cat))$ 
is a subcategory of
${\rm Op}_{\infty/\mathcal{O}^{\otimes}}^{\rm gen}\times_{\cat}\{S\}$
for any $S\in \cat$.

By unstraightening this natural transformation,
we obtain the following proposition.

\begin{proposition}\label{proposition:Lax-R-LAd-sub-Op}
The $\infty$-category 
$\mathsf{Mon}_{\mathcal{O}}^{\rm lax,R}(\cat)^{\rm LAd}$
is equivalent to a subcategory of 
${\rm Op}_{\infty/\mathcal{O}^{\otimes}}^{\rm gen}$.
The inclusion functor 
$\mathsf{Mon}_{\mathcal{O}}^{\rm lax,R}(\cat)^{\rm LAd}
\hookrightarrow
{\rm Op}_{\infty/\mathcal{O}^{\otimes}}^{\rm gen}$
is a map of Cartesian fibrations
over $\cat$.
\end{proposition}

\subsection{$\mathsf{Mon}_{\mathcal{O}}^{\rm oplax,L}(\cat)^{\rm RAd}$
and ${\rm Op}_{\infty/\mathcal{O}^{\otimes,{\rm op}}}^{\rm gen,\vee}$}
\label{subsection:mon-RAd-opposite-op}

In this subsection
we introduce an $\infty$-category
${\rm Op}_{\infty}^{{\rm gen},\vee}$
which is equivalent to
${\rm Op}_{\infty}^{{\rm gen}}$
and 
show that 
$\mathsf{Mon}_{\mathcal{O}}^{\rm oplax,L}(\cat)^{\rm RAd}$
is equivalent to a subcategory of 
${\rm Op}_{\infty/\mathcal{O}^{\otimes,{\rm op}}}^{{\rm gen},\vee}$.

We define an $\infty$-category 
\[ {\rm Op}_{\infty}^{{\rm gen},\vee} \]
to be a subcategory of $\cat{}_{/{\rm Fin}_*^{\rm op}}$
whose objects are maps $\mathcal{E}^{\otimes}\to {\rm Fin}_*^{\rm op}$
such that the opposite $\mathcal{E}^{\otimes,{\rm op}}\to
{\rm Fin}_*$ is a generalized $\infty$-operad
and whose morphisms are maps
$\mathcal{E}^{\otimes}\to \mathcal{F}^{\otimes}$
over ${\rm Fin}_*^{\rm op}$ such that
the opposite $\mathcal{E}^{\otimes,{\rm op}}\to
\mathcal{F}^{\otimes,{\rm op}}$ is a map
of generalized $\infty$-operads.
We have a functor
\[ (-)_{\langle 0\rangle}: 
   {\rm Op}_{\infty/\mathcal{O}^{\otimes,{\rm op}}}^{{\rm gen},\vee}\to
   \cat \]
which assigns the fiber 
$\mathcal{E}^{\otimes}_{\langle 0\rangle}$ at $\langle 0\rangle
\in \mathcal{O}^{\otimes,{\rm op}}$
to a map 
$\mathcal{E}^{\otimes}\to \mathcal{O}^{\otimes,{\rm op}}$
in ${\rm Op}_{\infty}^{{\rm gen},\vee}$.
Since ${\rm Op}_{\infty}^{\rm gen,\vee}$
is equivalent to ${\rm Op}_{\infty}^{\rm gen}$,
we see that the map
$(-)_{\langle 0\rangle}: 
{\rm Op}_{\infty/\mathcal{O}^{\otimes,{\rm op}}}^{{\rm gen},\vee}\to\cat$
is a Cartesian fibration
by Lemma~\ref{lemma:0-fiber-cartesian-1}.

By definition,
there is a natural functor
${\rm Fun}^{\rm RAd}(S,{\rm Mon}_{\mathcal{O}}^{\rm oplax,L}(\cat))
\to {\rm Fun}(S,{\rm Op}_{\infty/\mathcal{O}^{\otimes,{\rm op}}}^{\vee})$.
In the same way as in the case of 
$\mathsf{Mon}_{\mathcal{O}}^{\rm oplax,L}(\cat)^{\rm LAd}$,
we can identify ${\rm Fun}(S,
{\rm Op}_{\infty/\mathcal{O}^{\otimes,\rm op}}^{\vee})$
with a subcategory of 
${\rm Op}_{\infty/\mathcal{O}^{\otimes,\rm op}}^{{\rm gen},\vee}\times_{\cat}
\{S\}$.
Hence we obtain 
a natural transformation
\[ {\rm Fun}^{\rm RAd}
  (-,{\rm Mon}_{\mathcal{O}}^{\rm oplax,L}(\cat))\to
  {\rm Op}_{\infty/\mathcal{O}^{\otimes,\rm op}}^{{\rm gen},\vee}
  \times_{\cat}\{-\} \]
in which
${\rm Fun}^{\rm RAd}(S,{\rm Mon}_{\mathcal{O}}^{\rm oplax,L}(\cat))$ 
is a subcategory of
${\rm Op}_{\infty/\mathcal{O}^{\otimes,\rm op}}^{\rm gen,\vee}\times_{\cat}\{S\}$
for any $S\in \cat$.

By unstraightening this natural transformation,
we obtain the following proposition.

\begin{proposition}\label{prop:oplax-L-RAd-sub-Op}
The $\infty$-category 
$\mathsf{Mon}_{\mathcal{O}}^{\rm oplax,L}(\cat)^{\rm RAd}$
is equivalent to a subcategory of
the $\infty$-category
${\rm Op}_{\infty/\mathcal{O}^{\otimes,{\rm op}}}^{{\rm gen},\vee}$.
The inclusion functor
$\mathsf{Mon}_{\mathcal{O}}^{\rm oplax,L}(\cat)^{\rm RAd}
\hookrightarrow
{\rm Op}_{\infty/\mathcal{O}^{\otimes,{\rm op}}}^{{\rm gen},\vee}$
is a map of Cartesian fibrations
over $\cat$.
\end{proposition}

\section{Mixed fibrations
of operadic modules}
\label{section:mixed-fibration-operadic-module}

Let $\mathcal{O}^{\otimes}$
be a coherent $\infty$-operad and
let $q: \mathcal{C}^{\otimes}\to \mathcal{O}$
be a map of $\infty$-operads.
In \S\ref{subsection:infinity-cat-operadic-modules}
we recalled the construction of 
the map 
$(\Phi,\sigma):
{\rm Mod}^{\mathcal{O}}(\mathcal{C})^{\otimes}
\to {\rm Alg}_{/\mathcal{O}}(\mathcal{C})\times
\mathcal{O}^{\otimes}$ of generalized $\infty$-operads.
If $\mathcal{C}^{\otimes}$ has a sufficient supply of colimits,
then this encodes a structure consisting of
$\mathcal{O}$-monoidal $\infty$-categories
${\rm Mod}_A^{\mathcal{O}}(\mathcal{C})^{\otimes}$
and lax $\mathcal{O}$-monoidal functors
$f^*: {\rm Mod}_B^{\mathcal{O}}(\mathcal{C})^{\otimes}\to
{\rm Mod}_A^{\mathcal{O}}(\mathcal{C})^{\otimes}$.
In \S\ref{subsection:construction-psi-tau}
we construct a map
$(\Psi,\tau): {\rm Mod}^{\mathcal{O}}(\mathcal{C})^{\otimes,\vee}
\to {\rm Alg}_{/\mathcal{O}}(\mathcal{C})\times
\mathcal{O}^{\otimes,{\rm op}}$
which encodes a structure consisting of
$\mathcal{O}$-monoidal $\infty$-categories
${\rm Mod}_A^{\mathcal{O}}(\mathcal{C})^{\otimes}$
and oplax $\mathcal{O}$-monoidal functors
$f_!: {\rm Mod}_B^{\mathcal{O}}(\mathcal{C})^{\otimes}\to
{\rm Mod}_A^{\mathcal{O}}(\mathcal{C})^{\otimes}$.
In \S\ref{subsection:functoriality-psi-tau}
we study a functoriality
of the construction of $(\Psi,\tau)$
for $\mathcal{C}^{\otimes}$.
We construct a coCartesian fibration
$\Theta: {\rm Mod}^{\mathcal{O},{\rm Triple}}
    (\catkappa)^{\otimes,\vee}
    \to
    {\rm Alg}_{\mathcal{O}}^{\rm Pair}
    (\catkappa)$
which is associated to
a functor given by
$(\mathcal{C}^{\otimes},A)\mapsto 
{\rm Mod}_A^{\mathcal{O}}(\mathcal{C})^{\otimes,\vee}$.

\subsection{Construction
of the map $(\Psi,\tau)$}
\label{subsection:construction-psi-tau}

Let $\kappa$ be an uncountable regular cardinal.
We assume that
$\mathcal{O}^{\otimes}$ is an essentially $\kappa$-small
coherent $\infty$-operad and
that $\mathcal{C}^{\otimes}$ is an $\mathcal{O}$-monoidal
$\infty$-category which is compatible with 
$\kappa$-small colimits.
In this subsection
we construct a map
$(\Psi,\tau): {\rm Mod}^{\mathcal{O}}(\mathcal{C})^{\otimes,\vee}
\to {\rm Alg}_{/\mathcal{O}}(\mathcal{C})\times
\mathcal{O}^{\otimes,{\rm op}}$
which encodes a structure consisting of
$\mathcal{O}$-monoidal $\infty$-categories
${\rm Mod}_A^{\mathcal{O}}(\mathcal{C})^{\otimes}$
and oplax $\mathcal{O}$-monoidal functors
$f_!: {\rm Mod}_B^{\mathcal{O}}(\mathcal{C})^{\otimes}\to
{\rm Mod}_A^{\mathcal{O}}(\mathcal{C})^{\otimes}$.

Recall that 
$\Phi: {\rm Mod}^{\mathcal{O}}(\mathcal{C})^{\otimes}
\to{\rm Alg}_{/\mathcal{O}}(\mathcal{C})$
is a Cartesian fibration.
By the straightening functor,
there is a functor
\[ {\rm St}(\Phi):
   \algc{/\mathcal{O}}^{\mathrm{op}}\longrightarrow
   \cat, \]
which associates to $A\in\algc{/\mathcal{O}}$
the $\infty$-category
${\rm Mod}_A^{\mathcal{O}}(\mathcal{C})^{\otimes}$.
We can lift this functor to a functor
\[ {\rm St}(\Phi):
   \algc{/\mathcal{O}}^{\mathrm{op}}\longrightarrow
    \mathrm{Mon}_{\mathcal{O}}^{\rm lax,R}(\cat).\]

In \cite{Haugseng2, HLN, Torii4},
we have proved that there is an equivalence
\[ \mathrm{Mon}_{\mathcal{O}}^{\mathrm{oplax}, \rm L}(\cat)
   \simeq 
   (\mathrm{Mon}_{\mathcal{O}}^{\mathrm{lax},\rm R}(\cat))^{\mathrm{op}} \]
of $\infty$-categories.
By using this equivalence,
we obtain a functor
\[ \algc{/\mathcal{O}}\longrightarrow
   (\mathrm{Mod}_{\mathcal{O}}^{\rm lax,R}
   (\cat))^{\mathrm{op}}
   \simeq
   \mathrm{Mon}_{\mathcal{O}}^{\rm oplax ,L}(\cat)
   \hookrightarrow 
   \mathrm{Mon}_{\mathcal{O}}^{\mathrm{oplax}}(\cat).\]

Now,
we recall a description of
${\rm Mon}_{\mathcal{O}}^{\rm oplax}(\cat)$.
Let ${\rm Op}_{\infty}^{\vee}$
be a full subcategory of ${\rm Op}_{\infty}^{\rm gen,\vee}$
spanned by those maps
$\mathcal{E}^{\otimes}\to{\rm Fin}_*^{\rm op}$
such that the opposite $\mathcal{E}^{\otimes,\rm op}\to
{\rm Fin}_*$ is an $\infty$-operad.
We regard ${\rm Mon}^{\rm oplax}_{\mathcal{O}}(\cat)$
as a full subcategory of 
${\rm Op}_{\infty/\mathcal{O}^{\otimes,{\rm op}}}^{\vee}$
spanned by those maps
$\mathcal{E}^{\otimes}\to \mathcal{O}^{\otimes,{\rm op}}$ 
which is a Cartesian fibration.
We also say that an object
of ${\rm Mon}^{\rm oplax}_{\mathcal{O}}(\cat)$
is an $\mathcal{O}$-monoidal $\infty$-category
and a morphism is an oplax $\mathcal{O}$-monoidal
functor.

For a coCartesian fibration $\mathcal{X}\to S$
which is classified by a functor $S\to \cat$,
we denote by $\mathcal{X}^{\vee}\to S^{\rm op}$
a Cartesian fibration which is 
classified by the same functor.
For an $\mathcal{O}$-monoidal
$\infty$-category
$\mathcal{E}^{\otimes}\to\mathcal{O}^{\otimes}$,
the Cartesian fibration
$\mathcal{E}^{\otimes,\vee}\to\mathcal{O}^{\otimes,{\rm op}}$
is an object 
in $\mathrm{Mon}_{\mathcal{O}}^{\mathrm{oplax}}(\cat)$.
Note that the underlying $\infty$-category
$\mathcal{E}^{\otimes,\vee}_X$ 
of $\mathcal{E}^{\otimes,\vee}\to\mathcal{O}^{\otimes,{\rm op}}$
is equivalent to the underlying $\infty$-category
$\mathcal{E}^{\otimes}_X$ of 
$\mathcal{E}^{\otimes}\to\mathcal{O}^{\otimes}$
for each $X\in\mathcal{O}$.

The functor
$\algc{/\mathcal{O}}\to
\mathrm{Mon}_{\mathcal{O}}^{\mathrm{oplax}}(\cat)$
induces a commutative diagram
\begin{align}\label{align:triangle-Psi-tau}
   \xymatrix{
     \mathrm{Mod}^{\mathcal{O}}(\mathcal{C})^{\otimes, \vee}
     \ar[rr]^{(\Psi,\tau)}\ar[dr]_{\Psi} &&
     \mathrm{Alg}_{/\mathcal{O}}(\mathcal{C})\times
     \mathcal{O}^{\otimes,{\rm op}} \ar[dl]^{\pi}\\
     & \mathrm{Alg}_{/\mathcal{O}}(\mathcal{C}),\\ 
}\end{align}
where 
$\Psi$
is a coCartesian fibration,
the map $\tau$ carries $\Psi$-coCartesian morphisms
to equivalences,
and $\pi$ is the projection.
For each $A\in {\rm Alg}_{/\mathcal{O}}(\mathcal{C})$,
the fiber of $(\Psi,\tau)$ at $A$ determines 
a Cartesian fibration
\[ \tau_A: {\rm Mod}_A^{\mathcal{O}}
   (\mathcal{C})^{\otimes,\vee}
   \to \mathcal{O}^{\otimes,{\rm op}}, \]
which is an $\mathcal{O}$-monoidal
$\infty$-category equivalent
to $\sigma_A: {\rm Mod}_A^{\mathcal{O}}(\mathcal{C})^{\otimes}
\to \mathcal{O}^{\otimes}$.
A morphism $f: A\to B$ in
${\rm Alg}_{/\mathcal{O}}(\mathcal{C})$
induces an oplax $\mathcal{O}$-monoidal functor
\[ f_!: {\rm Mod}_A^{\mathcal{O}}
   (\mathcal{C})^{\otimes,\vee}
   \to {\rm Mod}_B^{\mathcal{O}}
   (\mathcal{C})^{\otimes,\vee},\]
where the restriction $f_{!X}$
is a left adjoint to $f^*_X$
for each $X\in\mathcal{O}$.

In \cite[Definition~3.15]{Torii1}
we introduced the notion of mixed fibrations.
A map $g: X\to S\times T$ of $\infty$-categories
is a mixed fibration over $(S,T)$
if it satisfies the following conditions:

\begin{itemize}

\item
The map $g_S: X\to S$ is a coCartesian fibration
and $g$ preserves coCartesian morphisms,
where $g_S$ is the composite of $g$ with the projection
$S\times T\to S$.

\item
The map $g_T: X\to T$ is a Cartesian fibration
and $g$ preserves Cartesian morphisms,
where $g_T$ is the composite of $g$ with the projection
$S\times T\to T$.

\end{itemize}

By \cite[Proposition~3.25]{Torii1},
we obtain the following proposition.

\begin{proposition}
There is a mixed fibration
\[ (\Psi,\tau): \mathrm{Mod}^{\mathcal{O}}(\mathcal{C})^{\otimes,\vee}
         \longrightarrow
         \mathrm{Alg}_{/\mathcal{O}}(\mathcal{C})
         \times\mathcal{O}^{\otimes,{\rm op}}\]
over 
$({\rm Alg}_{/\mathcal{O}}(\mathcal{C}),
\mathcal{O}^{\otimes,{\rm op}})$.
For $A\in {\rm Alg}_{/\mathcal{O}}(\mathcal{C})$,
the fiber of $(\Psi,\tau)$ at $A$ 
is a Cartesian fibration
$\tau_A: {\rm Mod}_A^{\mathcal{O}}
   (\mathcal{C})^{\otimes,\vee}
   \to \mathcal{O}^{\otimes,{\rm op}}$,
which is an $\mathcal{O}$-monoidal $\infty$-category
equivalent to $\sigma_A$.
For $f: A\to B$ in
${\rm Alg}_{/\mathcal{O}}(\mathcal{C})$,
the induced functor
$f_!: {\rm Mod}_A^{\mathcal{O}}
   (\mathcal{C})^{\otimes,\vee}
   \to {\rm Mod}_B^{\mathcal{O}}
   (\mathcal{C})^{\otimes,\vee}$
is an oplax $\mathcal{O}$-monoidal functor,
in which the restriction $f_{!X}$
is a left adjoint to $f^*_X$
for each $X\in\mathcal{O}$.
\end{proposition}

\subsection{A functoriality
of the construction of $(\Psi,\tau)$}
\label{subsection:functoriality-psi-tau}

In this subsection
we consider a functoriality
of the construction of 
the map $(\Psi,\tau)$
for $\mathcal{C}^{\otimes}$. 
We construct a coCartesian fibration
$\Theta: {\rm Mod}^{\mathcal{O},{\rm Triple}}
    (\catkappa)^{\otimes,\vee}
    \to
    {\rm Alg}_{\mathcal{O}}^{\rm Pair}
    (\catkappa)$
which is associated to a functor given by
$(\mathcal{C}^{\otimes},A)\mapsto
{\rm Mod}_A^{\mathcal{O}}(\mathcal{C})^{\otimes,\vee}$.

Let 
\[ \catkappa \]
be the subcategory of $\cat$ 
spanned by those small $\infty$-categories which have
$\kappa$-small colimits, and
those functors which preserve $\kappa$-colimits
(\cite[Definition~4.8.1.1]{Lurie2}).
By \cite[Corollary~4.8.1.4]{Lurie2},
the $\infty$-category $\catkappa$
inherits a symmetric monoidal structure
from $\cat$,
that is, 
there is a coCartesian fibration
\[ \catkappa{}^{\boxtimes} \longrightarrow {\rm Fin}_* \]
of $\infty$-operads
such that the inclusion functor
$\catkappa{}^{\boxtimes}\hookrightarrow
\cat^{\times}$ is a map
of $\infty$-operads.
We have an $\infty$-category ${\rm Alg}_{\mathcal{O}}
(\catkappa)$ of $\mathcal{O}$-algebra objects 
in $\catkappa$.
The objects 
are identified
with $\mathcal{O}$-monoidal $\infty$-categories
which are compatible with $\kappa$-small colimits
by \cite[Remark~4.8.1.9]{Lurie2}.

We have a functor
\[ {\rm Alg}_{/\mathcal{O}}(-):
   {\rm Alg}_{\mathcal{O}}(\catkappa)\longrightarrow
   \cat \]
which associates to $\mathcal{C}^{\otimes}$
the $\infty$-category
${\rm Alg}_{/\mathcal{O}}(\mathcal{C})$.
By unstraightening,
we obtain a coCartesian fibration
\[ \upsilon:
   {\rm Alg}_{\mathcal{O}}^{\rm Pair}(\catkappa)
   \longrightarrow
   {\rm Alg}_{\mathcal{O}}(\catkappa), \]
where the objects of ${\rm Alg}_{\mathcal{O}}^{\rm Pair}
(\catkappa)$
are pairs $(\mathcal{C}^{\otimes},A)$
of an $\mathcal{O}$-monoidal $\infty$-category $\mathcal{C}^{\otimes}$
that is compatible with $\kappa$-small colimits
and $A$ is an $\mathcal{O}$-algebra object
in $\mathcal{C}^{\otimes}$. 

We have a functor
\[ {\rm Mod}^{\mathcal{O}}(-)^{\otimes,\vee}:
   {\rm Alg}_{\mathcal{O}}(\catkappa)
   \longrightarrow
   \cat \]
which associates to $\mathcal{C}^{\otimes}$
the $\infty$-category
${\rm Mod}^{\mathcal{O}}(\mathcal{C})^{\otimes,\vee}$.
We denote by
\[ \omega: 
    {\rm Mod}^{\mathcal{O},{\rm Triple}}
    (\catkappa)^{\otimes,\vee}
    \longrightarrow
    {\rm Alg}_{\mathcal{O}}(\catkappa) \]
the associated coCartesian fibration
by unstraightening.
The objects of 
${\rm Mod}^{\mathcal{O},{\rm Triple}}
(\catkappa)^{\otimes,\vee}$
are triples $(\mathcal{C}^{\otimes},A,M)$ of
an $\mathcal{O}$-monoidal $\infty$-category
$\mathcal{C}^{\otimes}$ that is compatible
with $\kappa$-small colimits,
an $\mathcal{O}$-algebra object $A$,
and an $\mathcal{O}$-$A$-module $M$
in $\mathcal{C}^{\otimes}$. 

The functor
${\rm Mod}^{\mathcal{O}}(-)^{\otimes,\vee}$
lifts to a functor
\[ {\rm Mod}^{\mathcal{O}}(-)^{\otimes,\vee}:
   {\rm Alg}_{\mathcal{O}}(\catkappa)
   \longrightarrow
   {\rm Op}_{\infty/\mathcal{O}^{\otimes,{\rm op}}}^{{\rm gen},\vee}, \]
which fits into the following commutative diagram
\[ \xymatrix{
     & {\rm Op}_{\infty/\mathcal{O}^{\otimes,\rm op}}^{\rm gen,\vee}
       \ar[d]^{(-)_{\langle 0\rangle}}\\
     {\rm Alg}_{\mathcal{O}}(\catkappa)
      \ar[ur]^{{\rm Mod}^{\mathcal{O}}(-)^{\otimes,\vee}}
      \ar[r]_-{{\rm Alg}_{/\mathcal{O}}(-)}&
      \cat.\\
   }\]   
By unstraightening,
this induces the following commutative diagram
\begin{align}\label{align:non-multiplicative-fundamental-diagram}
 \xymatrix{
    {\rm Mod}^{\mathcal{O},{\rm Triple}}
    (\catkappa)^{\otimes,\vee}
    \ar[rr]^{(\Theta,\rho)}
    \ar[dr]^{\Theta}\ar@/_20pt/[ddr]_{\omega}&&
    {\rm Alg}_{\mathcal{O}}^{\rm Pair}
    (\catkappa)
    \times \mathcal{O}^{\otimes,{\rm op}}
    \ar[dl]_{\pi}\ar@/^20pt/[ddl]^{\upsilon\circ\pi}\\
    &{\rm Alg}_{\mathcal{O}}^{\rm Pair}
     (\catkappa)\ar[d]_{\upsilon}&\\
    &{\rm Alg}_{\mathcal{O}}(\catkappa),&
   }
\end{align}
where 
the map
$\Theta$
carries $\omega$-coCartesian morphisms to
$\upsilon$-coCartesian morphisms,
the map
$\rho$
carries $\omega$-coCartesian morphisms to equivalences,
and $\pi$ is the projection.

We would like to show that
$\Theta$ is a coCartesian fibration.
For this purpose,
we consider the following situation.
Suppose that we have a commutative diagram
\[ \xymatrix{
    \mathcal{X}\ar[rr]^{f}\ar[dr]_p&&
    \mathcal{Y}\ar[dl]^q\\
    & \mathcal{Z}& \\
}\]
of $\infty$-categories.
We assume that $p$ and $q$ are coCartesian fibrations,
that $f_z: \mathcal{X}_z\to \mathcal{Y}_z$
is a coCartesian fibration for each $z\in \mathcal{Z}$,
and that $f$ carries $p$-coCartesian morphisms
to $q$-coCartesian morphisms.
In this situation
we obtain a commutative diagram
\[ \begin{array}{ccc}
     \mathcal{X}_z 
     & \stackrel{\phi_!}{\longrightarrow} &
     \mathcal{X}_{z'} \\
     \mbox{$\scriptstyle f_z$}
     \bigg\downarrow
     \phantom{\mbox{$\scriptstyle f_z$}} 
     & & 
     \phantom{\mbox{$\scriptstyle f_{z'}$}}
     \bigg\downarrow
     \mbox{$\scriptstyle f_{z'}$} \\
     \mathcal{Y}_z 
     & \stackrel{\phi_!}{\longrightarrow} &
     \mathcal{Y}_{z'} \\
   \end{array}\]
for any morphism $\phi: z\to z'$ in $\mathcal{Z}$.

The dual form of the following lemma was proved
in \cite[Proposition~8.3]{GHN}. 

\begin{lemma}[{\cite[Proposition~8.3]{GHN}}]
\label{lemma:sufficient-coCartesian-general}
If $\phi_!: \mathcal{X}_z\to\mathcal{X}_{z'}$ 
carries $f_z$-coCartesian morphisms
to $f_{z'}$-coCartesian morphisms
for any morphism $\phi: z\to z'$ in $\mathcal{Z}$,
then $f$ is a coCartesian fibration.
\end{lemma}

\begin{proposition}
\label{prop:Theta-coCartesian}
The functor
$\Theta$ is a coCartesian fibration.
\end{proposition}

\proof
For any $\mathcal{C}^{\otimes}\in
{\rm Alg}_{\mathcal{O}}(\catkappa)$,
the map $\Theta_{\mathcal{C}}:
{\rm Mod}^{\mathcal{O},\rm Triple}(\catkappa)^{\otimes,\vee}_{\mathcal{C}}
\to {\rm Alg}_{\mathcal{O}}(\catkappa)_{\mathcal{C}}$
is identified with
the map $\Psi:
{\rm Mod}^{\mathcal{O}}(\mathcal{C})^{\otimes,\vee}\to
{\rm Alg}_{/\mathcal{O}}(\mathcal{C})$
which is a coCartesian fibration.
For a strong $\mathcal{O}$-monoidal
functor $F: \mathcal{C}^{\otimes}\to\mathcal{D}^{\otimes}$,
we have a commutative diagram
\[ \begin{array}{ccc}
   {\rm Mod}^{\mathcal{O}}(\mathcal{C})^{\otimes,\vee}
   &\stackrel{F_*}{\longrightarrow}&
   {\rm Mod}^{\mathcal{O}}(\mathcal{D})^{\otimes,\vee}   \\
   \mbox{$\scriptstyle \Theta_{\mathcal{C}}$}
   \bigg\downarrow
   \phantom{\mbox{$\scriptstyle \Theta_{\mathcal{C}}$}}
    & & 
   \phantom{\mbox{$\scriptstyle \Theta_{\mathcal{D}}$}}
   \bigg\downarrow
   \mbox{$\scriptstyle \Theta_{\mathcal{D}}$} \\
   {\rm Alg}_{/\mathcal{O}}(\mathcal{C})
   &\stackrel{}{\longrightarrow}&
   {\rm Alg}_{/\mathcal{O}}(\mathcal{D}).\\
   \end{array}\]
By Lemma~\ref{lemma:sufficient-coCartesian-general},
in order to prove the proposition,
it suffices to show that 
$F_*$ carries $\Theta_{\mathcal{C}}$-coCartesian morphisms
to $\Theta_{\mathcal{D}}$-coCartesian morphisms.
This follows from
the following commutative diagram
\[ \begin{array}{ccc}
   {\rm Mod}_A^{\mathcal{O}}(\mathcal{C})^{\otimes,\vee}
   &\stackrel{F_*}{\longrightarrow}&
   {\rm Mod}^{\mathcal{O}}(\mathcal{D})^{\otimes,\vee}   \\
   \mbox{$\scriptstyle f_!$}
   \bigg\downarrow
   \phantom{\mbox{$\scriptstyle f_!$}} 
   & & 
   \phantom{\mbox{$\scriptstyle F(f)_!$}}
   \bigg\downarrow
   \mbox{$\scriptstyle F(f)_!$}  \\
   {\rm Mod}_B^{\mathcal{O}}(\mathcal{C})^{\otimes,\vee}
   &\stackrel{F_*}{\longrightarrow}&
   {\rm Mod}_{F(B)}^{\mathcal{O}}(\mathcal{C})^{\otimes,\vee}\\
   \end{array}\]
for any $f: A\to B$ in ${\rm Alg}_{/\mathcal{O}}(\mathcal{C})$.
\qed

\section{Duoidal $\infty$-categories of operadic modules}
\label{section:duoidal-operadic-module}

In this section
we will construct duoidal structures
on $\infty$-categories of operadic modules.
We fix an uncountable regular cardinal $\kappa$
and an essentially $\kappa$-small coherent $\infty$-operad
$\mathcal{O}^{\otimes}$
throughout this section
unless otherwise stated.
Let $\mathcal{P}^{\otimes}$ be an $\infty$-operad. 
We take a $\mathcal{P}\otimes\mathcal{O}$-monoidal 
$\infty$-category
$\mathcal{C}^{\otimes}$
that is compatible with $\kappa$-small colimits,
and a $\mathcal{P}\otimes\mathcal{O}$-algebra 
object $A$ in $\mathcal{C}^{\otimes}$.
We shall show that
the $\infty$-category ${\rm Mod}_A^{\mathcal{O}}(\mathcal{C})$
of $\mathcal{O}$-$A$-modules in $\mathcal{C}^{\otimes}$
has a structure of a $(\mathcal{P},\mathcal{O})$-duoidal
$\infty$-category
(Theorem~\ref{thm:main-functor}). 

\subsection{Mixed $(\mathcal{P},\mathcal{Q})$-monoidal
$\infty$-categories}

Let $\mathcal{P}^{\otimes}$ and $\mathcal{Q}^{\otimes}$
be $\infty$-operads.
In this subsection
we recall a formulation 
of $(\mathcal{P},\mathcal{Q})$-duoidal $\infty$-categories
as mixed $(\mathcal{P},\mathcal{Q})$-monoidal
$\infty$-categories.

We defined a $(\mathcal{P},\mathcal{Q})$-duoidal
$\infty$-category as a $\mathcal{P}$-monoid
object in the Cartesian symmetric monoidal
$\infty$-category ${\rm Mon}_{\mathcal{Q}}^{\rm oplax}(\cat)$.
There are three formulations
of the $\infty$-category of 
$(\mathcal{P},\mathcal{Q})$-duoidal $\infty$-categories 
depending on which kinds of functors we choose.
In this paper it is convenient
to formulate $(\mathcal{P},\mathcal{Q})$-duoidal $\infty$-categories
as mixed $(\mathcal{P},\mathcal{Q})$-monoidal $\infty$-categories.

\begin{definition}[{cf.~\cite[Definition~3.11]{Torii3}}]
\rm
A mixed $(\mathcal{P},\mathcal{Q})$-monoidal
$\infty$-category 
is a mixed fibration
\[ \mathcal{D}^{\otimes}\longrightarrow 
      \mathcal{P}^{\otimes}\times\mathcal{Q}^{\otimes,\rm op} \]
over $(\mathcal{P}^{\otimes},\mathcal{Q}^{\otimes,\rm op})$
which satisfies the following conditions:

\begin{itemize}

\item
For each $p\simeq p_1\oplus\cdots\oplus p_m
\in\mathcal{P}^{\otimes}$,
the Segal morphism
\[ \mathcal{D}^{\otimes}_p\longrightarrow
   \mathcal{D}^{\otimes}_{p_1}\times_{\mathcal{Q}^{\otimes,\rm op}}
   \cdots\times_{\mathcal{Q}^{\otimes,\rm op}}
   \mathcal{D}^{\otimes}_{p_m} \]
is an equivalence.

\item
For each $q\simeq q_1\oplus\cdots\oplus q_n
\in\mathcal{Q}^{\otimes}$,
the Segal morphism
\[ \mathcal{D}^{\otimes}_q\longrightarrow
   \mathcal{D}^{\otimes}_{q_1}\times_{\mathcal{P}^{\otimes}}
   \cdots\times_{\mathcal{P}^{\otimes}}
   \mathcal{D}^{\otimes}_{q_n} \]
is an equivalence.

\end{itemize}

\end{definition}


We can define a bilax $(\mathcal{P},\mathcal{Q})$-monoidal
functor between mixed $(\mathcal{P},\mathcal{Q})$-monoidal
$\infty$-categories.
By \cite[Theorem~3.12]{Torii3},
the $\infty$-category
${\rm Mon}_{\mathcal{P}}^{\rm lax}({\rm Mon}_{Q}^{\rm oplax}(\cat))$
is equivalent to
the $\infty$-category of mixed $(\mathcal{P},\mathcal{Q})$-monoidal
$\infty$-categories and bilax $(\mathcal{P},\mathcal{Q})$-monoidal
functors.
In this paper
we identify a $(\mathcal{P},\mathcal{Q})$-duoidal
$\infty$-category with the corresponding
mixed $(\mathcal{P},\mathcal{Q})$-monoidal $\infty$-category.

\subsection{Construction of a functor $\Theta^{\boxtimes}$}
\label{subsection:construction-psi-boxtimes}

In this subsection we construct a map
$\Theta^{\boxtimes}$ of generalized $\infty$-operads which encodes
a multiplicative structure on the triples 
$(\mathcal{C}^{\otimes},A,M)$,
where $\mathcal{C}^{\otimes}$ is an $\mathcal{O}$-monoidal
$\infty$-category that is compatible with $\kappa$-small colimits,
$A$ is an $\mathcal{O}$-algebra object,
and $M$ is an $\mathcal{O}$-$A$-module object in $\mathcal{C}^{\otimes}$.
In the next subsection
we will show that $\Theta^{\boxtimes}$ is a coCartesian fibration 
of generalized $\infty$-operads.

First,
we shall construct a coCartesian fibration
$\upsilon^{\boxtimes}:
   {\rm Alg}_{\mathcal{O}}^{\rm Pair}(\catkappa)^{\boxtimes}
   \to
   {\rm Alg}_{\mathcal{O}}(\catkappa)^{\boxtimes}$
of $\infty$-operads,
which is an extension of 
the map
$\upsilon:
   {\rm Alg}_{\mathcal{O}}^{\rm Pair}(\catkappa)
   \to
   {\rm Alg}_{\mathcal{O}}(\catkappa)$.

The $\infty$-category ${\rm Alg}_{\mathcal{O}}
(\catkappa)$ has a symmetric monoidal structure
by pointwise multiplication
(\cite[Example~3.2.4.4]{Lurie2}).
We denote by
\[ {\rm Alg}_{\mathcal{O}}(\catkappa)^{\boxtimes}
   \longrightarrow {\rm Fin}_* \]
the associated coCartesian fibration of $\infty$-operads.
The functor
${\rm Alg}_{/\mathcal{O}}(-):
   {\rm Mon}_{\mathcal{O}}(\cat)
   \to \cat$
preserves finite products and
thus induces a functor
${\rm Alg}_{/\mathcal{O}}(-)^{\times}:
{\rm Mon}_{\mathcal{O}}(\cat)^{\times}\to\cat$,
which is a ${\rm Mon}_{\mathcal{O}}(\cat)$-monoid
object in $\cat$.
The map
$\catkappa{}^{\boxtimes}\to\cat^{\times}$
of $\infty$-operads
induces a lax symmetric monoidal functor
${\rm Alg}_{\mathcal{O}}(\catkappa)^{\boxtimes}
\to {\rm Alg}_{\mathcal{O}}(\cat)^{\times}\simeq
{\rm Mon}_{\mathcal{O}}(\cat)^{\times}$.
By composing these two functors,
we obtain a functor 
${\rm Alg}_{/\mathcal{O}}(-)^{\boxtimes}:
{\rm Alg}_{\mathcal{O}}(\catkappa)^{\boxtimes}
   \to \cat$,
which is an ${\rm Alg}_{\mathcal{O}}(\catkappa)$-monoid
object in $\cat$.
We define a map
\[ \upsilon^{\boxtimes}: {\rm Alg}_{\mathcal{O}}^{\rm Pair}
   (\catkappa)^{\boxtimes}
   \longrightarrow 
   {\rm Alg}_{/\mathcal{O}}
   (\catkappa)^{\boxtimes} \]
of $\infty$-operads
to be the associated coCartesian fibration
by unstraightening.

Next,
we shall construct a functor
${\rm Mod}^{\mathcal{O}}(-)^{\boxtimes,\otimes,\vee}:
   {\rm Alg}_{\mathcal{O}}(\catkappa)^{\boxtimes}
    \to
   {\rm Op}_{\infty/\mathcal{O}^{\otimes,\rm op}}^{\rm gen,\vee}$,
which is an extension of
the functor 
${\rm Mod}^{\mathcal{O}}(-)^{\otimes,\vee}:
   {\rm Alg}_{\mathcal{O}}(\catkappa)
    \to
   {\rm Op}_{\infty/\mathcal{O}^{\otimes,\rm op}}^{\rm gen,\vee}$,
and show that it is an
${\rm Alg}_{\mathcal{O}}(\catkappa)$-monoid object
in ${\rm Op}_{\infty/\mathcal{O}^{\otimes,\rm op}}^{\rm gen,\vee}$.

We have a functor
${\rm Mod}^{\mathcal{O}}(-)^{\otimes}:
   {\rm Mon}_{\mathcal{O}}(\cat)\to
   {\rm Op}_{\infty/\mathcal{O}^{\otimes}}^{\rm gen}$
which associates to $\mathcal{C}^{\otimes}$
the map 
$\sigma: {\rm Mod}^{\mathcal{O}}(\mathcal{C})^{\otimes}\to
\mathcal{O}^{\otimes}$
of generalized $\infty$-operads.
The functor ${\rm Mod}^{\mathcal{O}}(-)^{\otimes}$
preserves finite products and hence it extends to
a functor
${\rm Mod}^{\mathcal{O}}(-)^{\times,\otimes}:
    {\rm Mon}_{\mathcal{O}}(\cat)^{\times}
    \to
    {\rm Op}_{\infty}^{\rm gen}{}_{/\mathcal{O}^{\otimes}}$,
which is a ${\rm Mon}_{\mathcal{O}}(\cat)$-monoid
object in ${\rm Op}_{\infty}^{\rm gen}{}_{/\mathcal{O}^{\otimes}}$.
By composing 
${\rm Mod}^{\mathcal{O}}(-)^{\times,\otimes}$
with the lax symmetric monoidal functor
${\rm Alg}_{\mathcal{O}}(\catkappa)^{\boxtimes}
   \to
   {\rm Alg}_{\mathcal{O}}(\cat)^{\times}
   \simeq
   {\rm Mon}_{\mathcal{O}}(\cat)^{\times}$,
we obtain a functor
\[ {\rm Mod}^{\mathcal{O}}(-)^{\boxtimes,\otimes}:
   {\rm Alg}_{\mathcal{O}}(\catkappa)^{\boxtimes}
    \longrightarrow
   {\rm Op}_{\infty}^{\rm gen}{}_{/\mathcal{O}^{\otimes}}, \]
which is an ${\rm Alg}_{\mathcal{O}}(\catkappa)$-monoid
object in ${\rm Op}_{\infty}^{\rm gen}{}_{/\mathcal{O}^{\otimes}}$.

The $\infty$-category
${\rm Mon}_{\mathcal{O}}^{\rm lax,R}(\cat)^{\rm LAd}$
is equivalent to a subcategory
of ${\rm Op}_{\infty/\mathcal{O}^{\otimes}}^{\rm gen}$
by Proposition~\ref{proposition:Lax-R-LAd-sub-Op}.
The functor
${\rm Mod}^{\mathcal{O}}(-)^{\boxtimes,\otimes}$
factors through the subcategory
${\rm Mon}_{\mathcal{O}}^{\rm lax,R}(\cat)^{\rm LAd}$.
By Propositions~\ref{prop:monoidal-adjointable-equivalence}
and 
\ref{prop:oplax-L-RAd-sub-Op},
we obtain a functor
\[ {\rm Mod}^{\mathcal{O}}(-)^{\boxtimes,\otimes,\vee}:
   {\rm Alg}_{\mathcal{O}}(\catkappa)^{\boxtimes}
    \longrightarrow
    {\rm Op}_{\infty}^{{\rm gen}, \vee}{}_{/\mathcal{O}^{\otimes,{\rm op}}}. \]

\begin{lemma}
The functor
${\rm Mod}^{\mathcal{O}}(-)^{\boxtimes,\otimes,\vee}$
is an ${\rm Alg}_{\mathcal{O}}(\catkappa)$-monoid object
of ${\rm Op}_{\infty}^{{\rm gen}, \vee}{}_{/\mathcal{O}^{\otimes,{\rm op}}}$.
\end{lemma}

\proof
Let $(\mathcal{C}_1^{\otimes},\ldots,\mathcal{C}_n^{\otimes})
\in {\rm Alg}_{\mathcal{O}}(\catkappa)^{\boxtimes}$.
By definition,
the image of $(\mathcal{C}_1^{\otimes},\ldots,\mathcal{C}_n^{\otimes})$
under the functor
${\rm Mod}^{\mathcal{O}}(-)^{\boxtimes,\otimes,\vee}:
{\rm Alg}_{\mathcal{O}}(\catkappa)^{\boxtimes}
    \to
    {\rm Op}_{\infty}^{{\rm gen}, \vee}{}_{/\mathcal{O}^{\otimes,{\rm op}}}$
is given by a Cartesian fibration
\[ {\rm Mod}^{\mathcal{O}}
   (\mathcal{C}_1\times_{\mathcal{O}}\cdots
    \times_{\mathcal{O}}\mathcal{C}_n)^{\otimes,\vee} 
   \longrightarrow \mathcal{O}^{\otimes,{\rm op}}\]
of generalized $\infty$-operads.
We can see that 
it is equivalent to
a product of the Cartesian fibrations
${\rm Mod}^{\mathcal{O}}(\mathcal{C}_i)^{\otimes,\vee}\to
\mathcal{O}^{\otimes,{\rm op}}$
in ${\rm Op}_{\infty}^{{\rm gen},\vee}{}_{/\mathcal{O}^{\otimes,{\rm op}}}$
for $1\le i\le n$. 
\qed

\bigskip

By unstraightening the functor
${\rm Mod}^{\mathcal{O}}(-)^{\boxtimes,\otimes,\vee}$,
we obtain the following commutative diagram
\[ \xymatrix{
    {\rm Mod}^{\mathcal{O},{\rm Triple}}
    (\catkappa)^{\boxtimes,\otimes,\vee}
    \ar[rr]^{(\Theta^{\boxtimes},\tau^{\boxtimes})}
    \ar[dr]^{\Theta^{\boxtimes}}\ar@/_20pt/[ddr]_{\omega^{\boxtimes}}&&
    {\rm Alg}_{\mathcal{O}}^{\rm Pair}
    (\catkappa)^{\boxtimes}
    \times \mathcal{O}^{\otimes,{\rm op}}
    \ar[dl]_{\pi}\ar@/^20pt/[ddl]^{\upsilon^{\boxtimes}\circ\pi}\\
    &{\rm Alg}_{\mathcal{O}}^{\rm Pair}
     (\catkappa)^{\boxtimes}\ar[d]_{\upsilon^{\boxtimes}}&\\
    &{\rm Alg}_{\mathcal{O}}(\catkappa)^{\boxtimes},&
   }\]
where the maps 
$\omega^{\boxtimes}$
and 
$\upsilon^{\boxtimes}$
are coCartesian fibrations of generalized $\infty$-operads,
the map
$\Theta^{\boxtimes}$
is a map of generalized $\infty$-operads which
carries $\omega^{\boxtimes}$-coCartesian morphisms to
$\upsilon^{\boxtimes}$-coCartesian morphisms,
and the map
$\tau^{\boxtimes}$
carries $\omega^{\boxtimes}$-coCartesian morphisms to equivalences.

\subsection{The coCartesian fibration $\Theta^{\boxtimes}$}
\label{subsection:cocartesian-fibration-psi-boxtimes}

In \S\ref{subsection:construction-psi-boxtimes}
we have constructed a map
\[ \Theta^{\boxtimes}: {\rm Mod}^{\mathcal{O,{\rm Triple}}}
   (\catkappa)^{\boxtimes,\otimes,\vee}
   \longrightarrow 
   {\rm Alg}_{\mathcal{O}}^{\rm pair}
   (\catkappa)^{\boxtimes}
\]
of generalized $\infty$-operads.
In this subsection we shall show that
$\Theta^{\boxtimes}$ is a coCartesian fibration.

For simplicity,
we set $\mathcal{M}^T={\rm Mod}^{\mathcal{O},{\rm Triple}}
(\catkappa)^{\boxtimes,\otimes,\vee}$,
$\mathcal{A}^P={\rm Alg}_{\mathcal{O}}^{\rm Pair}
(\catkappa)^{\boxtimes}$,
and
$\mathcal{A}={\rm Alg}_{\mathcal{O}}
(\catkappa)^{\boxtimes}$.
Recall that we have the commutative diagram
\[ \xymatrix{
    \mathcal{M}^T\ar[rr]^{\Theta^{\boxtimes}}\ar[dr]_{\omega^{\boxtimes}}&&
    \mathcal{A}^P\ar[dl]^{\upsilon^{\boxtimes}}\\
    & \mathcal{A},& \\
}\]
where $\omega^{\boxtimes}$ and 
$\upsilon^{\boxtimes}$ are coCartesian fibrations,
and $\Theta^{\boxtimes}$ carries $\omega^{\boxtimes}$-coCartesian morphisms
to $\upsilon^{\boxtimes}$-coCartesian morphisms.
Furthermore,
for any object 
$C\in\mathcal{A}$,
the induced functor
$\Theta_C^{\boxtimes}: \mathcal{M}^T_C\to \mathcal{A}^P_C$
is a coCartesian fibration.

Let $\phi: 
(\mathcal{C}_1^{\otimes},\ldots,\mathcal{C}_n^{\otimes})
\to \mathcal{D}^{\otimes}$ be an active morphism in 
$\mathcal{A}$,
where $\mathcal{C}_1^{\otimes},\ldots,\mathcal{C}_n^{\otimes},
\mathcal{D}^{\otimes}$ are objects in 
the underlying $\infty$-category
${\rm Alg}_{\mathcal{O}}(\catkappa)$.
Then we have a commutative diagram
\[ \begin{array}{ccc}
    {\rm Mod}^{\mathcal{O}}
    (\mathcal{C}_1\times_{\mathcal{O}}\cdots
     \times_{\mathcal{O}}\mathcal{C}_n)^{\otimes,\vee}
    &\stackrel{\overline{\phi}_*}{\longrightarrow}&
    {\rm Mod}^{\mathcal{O}}(\mathcal{D})^{\otimes,\vee}\\
    \mbox{$\scriptstyle p$}
    \bigg\downarrow
    \phantom{\mbox{$\scriptstyle p$}} 
    & & 
    \phantom{\mbox{$\scriptstyle q$}} 
    \bigg\downarrow
    \mbox{$\scriptstyle q$} \\
    {\rm Alg}_{/\mathcal{O}}
    (\mathcal{C}_1\times_{\mathcal{O}}\cdots
     \times_{\mathcal{O}}\mathcal{C}_n)
    &\stackrel{\underline{\phi}_*}{\longrightarrow} &
    {\rm Alg}_{/\mathcal{O}}(\mathcal{D}).\\
   \end{array}\]
In order to prove that 
$\Theta^{\boxtimes}$ is a coCartesian fibration,
it suffices to show that
$\overline{\phi}_*$ carries
$p$-coCartesian morphisms to $q$-coCartesian morphisms
for any $\phi$
by Lemma~\ref{lemma:sufficient-coCartesian-general}.

First,
we consider the case in which 
$\mathcal{D}=\mathcal{C}_1\boxtimes\cdots\boxtimes\mathcal{C}_n$
and 
$\phi: (\mathcal{C}_1,\ldots,\mathcal{C}_n)\to
\mathcal{C}_1\boxtimes\cdots\boxtimes\mathcal{C}_n$
is a multiplication map
in the symmetric monoidal
$\infty$-category ${\rm Alg}_{\mathcal{O}}(\catkappa)$.

\begin{lemma}
\label{lemma:f-shriek-boxproduct-compatibility}
When $\phi: (\mathcal{C}_1,\ldots,\mathcal{C}_n)\to
\mathcal{C}_1\boxtimes\cdots\boxtimes\mathcal{C}_n$
is a multiplication map,
the functor
$\overline{\phi}_*$
carries $p$-coCartesian morphisms
to $q$-coCartesian morphisms.
\end{lemma}

\proof
Let $f_i: A_i\to B_i$ be a map
in ${\rm Alg}_{/\mathcal{O}}(\mathcal{C}_i)$,
and let $M_i$ be an object
of ${\rm Mod}_{A_i}^{\mathcal{O}}(\mathcal{C}_i)_X^{\otimes,\vee}$
for $1\le i\le n$ and $X\in\mathcal{O}^{\otimes}$.
We have to show that the canonical map
\[ (f_1\boxtimes\cdots\boxtimes f_n)_!
   (M_1\boxtimes\cdots\boxtimes M_n)
   \stackrel{}{\longrightarrow}
   (f_{1!}M_1)\boxtimes\cdots\boxtimes
   (f_{n!}M_n)\]
is an equivalence
in ${\rm Mod}_{B_1\boxtimes\cdots\boxtimes B_n}^{\mathcal{O}}
(\mathcal{C}_1\boxtimes\cdots\boxtimes\mathcal{C}_n)^{\otimes,\vee}_X$.

In order to prove this,
it suffices to show the case
in which $f_2,\ldots,f_n$ are identities.
Thus, we will prove that
the canonical map
\[ (f_1\boxtimes{\rm id}\boxtimes\ldots\boxtimes{\rm id})_{!}
   (M_1\boxtimes M_2\boxtimes \ldots\boxtimes M_n)
   \longrightarrow
   (f_{1!}M_1)\boxtimes M_2\boxtimes \ldots \boxtimes M_n \]
is an equivalence.

For simplicity, we assume that $n=2$.
Since the functor 
$(f_1\boxtimes {\rm id})^*:
   {\rm Mod}_{B_1\boxtimes B_2}^{\mathcal{O}}
   (\mathcal{C}_1\boxtimes\mathcal{C}_2)^{\otimes,\vee}_X\to
   {\rm Mod}_{A_1\boxtimes A_2}^{\mathcal{O}}
   (\mathcal{C}_1\boxtimes\mathcal{C}_2)^{\otimes,\vee}_X$
is conservative,
it suffices to show that
the map
\[ 
    (f_1\boxtimes {\rm id})^*
   (f_1\boxtimes {\rm id})_{!}
   (M_1\boxtimes M_2)
   \stackrel{}{\longrightarrow}
   (f_1\boxtimes {\rm id})^*
    ((f_{1!}M_1)\boxtimes M_2)
\]
is an equivalence.
Note that there is an equivalence
$
     (f_1\boxtimes {\rm id})^*
   ((f_{1!}M_1)\boxtimes M_2)
   \simeq 
   (f_1^*f_{1!}M_1)\boxtimes M_2$.

By Proposition~\ref{prop:functor:alg-C-fun-I-C}
we have an equivalence
$f_1^*f_{1!}M_1\simeq
   {\rm colim}_{I(\mathcal{O})}
   D(B_1,M_1)$.
Since $\boxtimes$ preserves
small colimits separately in each variable,
we obtain an equivalence
$(f_1^*f_{1!}M_1)\boxtimes M_2
   \simeq
   {\rm colim}_{I(\mathcal{O})}
   (D(B_1,M_1)\boxtimes M_2)$.

On the other hand,
we have an equivalence
$
    (f_1\boxtimes {\rm id})^*
   (f_1\boxtimes {\rm id})_!
   (M_1\boxtimes M_2)
   \simeq 
   {\rm colim}_{I(\mathcal{O})}
   D(B_1\boxtimes A_2,
     M_1\boxtimes M_2)$
by Proposition~\ref{prop:functor:alg-C-fun-I-C}.
Since we have an equivalence
$D(B_1\boxtimes A_2,
     M_1\boxtimes M_2)
   \simeq
   D(B_1,M_1)\boxtimes M_2$
of functors,
we obtain the desired equivalence.
\qed

\bigskip

For any active morphism
$\phi: (\mathcal{C}_1^{\otimes},\ldots,\mathcal{C}_n^{\otimes})
\to \mathcal{D}^{\otimes}$
in ${\rm Alg}_{\mathcal{O}}(\catkappa)^{\boxtimes}$,
where $\mathcal{C}_1^{\otimes},\ldots,\mathcal{C}_n^{\otimes},
\mathcal{D}^{\otimes}$ are objects
in the underlying $\infty$-category ${\rm Alg}_{\mathcal{O}}
(\catkappa)$,
we can decompose $\phi$ as 
\[ 
   (\mathcal{C}_1^{\otimes},\cdots,\mathcal{C}_n^{\otimes})
   \stackrel{}{\longrightarrow} 
   \mathcal{C}_1^{\otimes}\boxtimes\cdots\boxtimes 
   \mathcal{C}_n^{\otimes}
   \stackrel{}{\longrightarrow}
   \mathcal{D}^{\otimes},\]      
where the first map is a multiplication map and 
the second map is a strong $\mathcal{O}$-monoidal
functor that preserves $\kappa$-small colimits.
By Lemma~\ref{lemma:f-shriek-boxproduct-compatibility}
and the proof of Proposition~\ref{prop:Theta-coCartesian}, 
we obtain the following proposition.

\begin{proposition}\label{prop:Phi-coCartesian-fibration}
The map
$\Theta^{\boxtimes}: {\rm Mod}^{\mathcal{O},{\rm Triple}}
(\catkappa)^{\boxtimes,\otimes,\vee}
\to {\rm Alg}_{/\mathcal{O}}^{\rm Pair}
(\catkappa)^{\boxtimes}$
is a coCartesian fibration
of generalized $\infty$-operads.
\end{proposition}

\subsection{Construction 
of duoidal $\infty$-categories
of operadic modules}
\label{subsection:main-results}

In this subsection
we prove the main theorem (Theorem~\ref{thm:main-functor}).
For any $\infty$-operad $\mathcal{P}^{\boxtimes}$,
any $\mathcal{P}\otimes\mathcal{O}$-monoidal
$\infty$-category
$\mathcal{C}^{\otimes}$ 
that is compatible
with $\kappa$-small colimits,
and any $\mathcal{P}\otimes\mathcal{O}$-algebra 
object $A$ in $\mathcal{C}^{\otimes}$, 
we show that 
the $\infty$-category ${\rm Mod}^{\mathcal{O}}_A(\mathcal{C})$
of $\mathcal{O}$-$A$-modules in $\mathcal{C}^{\otimes}$
has a structure of $(\mathcal{P},\mathcal{O})$-duoidal
$\infty$-category.

First, we consider the universal case
in which $\mathcal{P}^{\boxtimes}=
{\rm Alg}_{\mathcal{O}}^{\rm Pair}(\catkappa)^{\boxtimes}$.
We have a symmetric monoidal $\infty$-category
${\rm Alg}_{\mathcal{O}}^{\rm Pair}(\catkappa)^{\boxtimes}$
and regard it as an $\infty$-operad.
We have a functor 
\[ {\rm Mod}^{\mathcal{O},{\rm Triple}}
   (\catkappa)^{\otimes,\vee}:
   {\rm Alg}_{\mathcal{O}}^{\rm Pair}
   (\catkappa)
    \longrightarrow
   {\rm Mon}_{\mathcal{O}}^{\rm oplax}(\cat)\]
which assigns to a pair $(\mathcal{C}^{\otimes}, A)$
the $\mathcal{O}$-monoidal $\infty$-category 
${\rm Mod}_A^{\mathcal{O}}(\mathcal{C})^{\otimes,\vee}
\to \mathcal{O}^{\otimes,{\rm op}}$.
By Proposition~\ref{prop:Phi-coCartesian-fibration}
and the straightening functor,
we obtain the following theorem.

\begin{theorem}
Let $\kappa$ be an uncountable regular cardinal and
let $\mathcal{O}^{\otimes}$
be an essentially $\kappa$-small coherent $\infty$-operad. 
Then we have a mixed $({\rm Alg}_{\mathcal{O}}^{\rm pair}(\catkappa),
\mathcal{O})$-monoidal $\infty$-category
\[ {\rm Mod}^{\mathcal{O},{\rm Triple}}
   (\catkappa)^{\boxtimes,\otimes,\vee}
   \longrightarrow
   {\rm Alg}_{\mathcal{O}}^{\rm pair}(\catkappa)\times
   \mathcal{O}^{\otimes,\rm op}, \]
in which 
the underlying $\infty$-category
${\rm Mod}^{\mathcal{O},{\rm Triple}}
   (\catkappa)^{\boxtimes,\otimes,\vee}_{(\mathcal{C}^{\otimes},A,X)}$
is equivalent to
${\rm Mod}_A^{\mathcal{O}}(\mathcal{C})^{\otimes}_X$
for $(\mathcal{C}^{\otimes},A)\in 
{\rm Alg}_{\mathcal{O}}^{\rm Pair}(\catkappa)$
and $X\in \mathcal{O}$.
\end{theorem}

\proof
By Proposition~\ref{prop:Phi-coCartesian-fibration}
and the straightening functor,
we can extend 
${\rm Mod}^{\mathcal{O},{\rm Triple}}
(\catkappa)^{\otimes,\vee}$
to a functor
\[ {\rm Mod}^{\mathcal{O},{\rm Triple}}
   (\catkappa)^{\boxtimes,\otimes,\vee}:
   {\rm Alg}_{\mathcal{O}}^{\rm Pair}
   (\catkappa)^{\boxtimes}
    \longrightarrow
   {\rm Mon}^{\rm oplax}(\cat) \]
which is an ${\rm Alg}_{\mathcal{O}}^{\rm Pair}
(\catkappa)$-monoid
object in ${\rm Mon}^{\rm oplax}(\cat)$.
\qed

\bigskip

Next,
we consider the case in which 
$\mathcal{P}^{\boxtimes}$ is any $\infty$-operad.
We set
\[ {\rm Alg}_{\mathcal{P}\otimes\mathcal{O}}^{\rm Pair}
    (\catkappa)
   ={\rm Alg}_{\mathcal{P}}
    ({\rm Alg}_{\mathcal{O}}^{\rm Pair}
    (\catkappa)), \]
where the objects of
${\rm Alg}_{\mathcal{P}\otimes\mathcal{O}}^{\rm Pair}
    (\catkappa)$
are pairs $(\mathcal{C}^{\otimes},A)$
of a $\mathcal{P}\otimes\mathcal{O}$-monoidal 
$\infty$-category $\mathcal{C}^{\otimes}$
that is compatible with $\kappa$-small colimits
and a $\mathcal{P}\otimes\mathcal{O}$-algebra
object $A$ in $\mathcal{C}^{\otimes}$.

The functor
${\rm Mod}^{\mathcal{O},{\rm Triple}}
   (\catkappa)^{\boxtimes,\otimes,\vee}$
induces a lax symmetric monoidal
functor
\[ 
   {\rm Alg}_{\mathcal{O}}^{\rm Pair}
   (\catkappa)^{\boxtimes}
    \longrightarrow
   {\rm Mon}^{\rm oplax}(\cat)^{\times}. \]
By applying the functor ${\rm Alg}_{\mathcal{P}}(-)$
to this functor,
we obtain the following theorem.

\begin{theorem}\label{thm:main-functor}
Let $\kappa$ be an uncountable regular cardinal and
let $\mathcal{O}^{\otimes}$ be an essentially $\kappa$-small
coherent $\infty$-operad.
For any pair $(\mathcal{C}^{\otimes},A)$
of a $\mathcal{P}\otimes\mathcal{O}$-monoidal
$\infty$-category $\mathcal{C}^{\otimes}$
that is compatible with $\kappa$-small colimits
and a $\mathcal{P}\otimes\mathcal{O}$-algebra
object $A$ in $\mathcal{C}^{\otimes}$,
we have a mixed $(\mathcal{P},\mathcal{O})$-monoidal
$\infty$-category
\[ {\rm Mod}_A^{\mathcal{O}}(\mathcal{C})^{\boxtimes, \otimes,\vee}
   \longrightarrow \mathcal{P}^{\boxtimes}\times
   \mathcal{O}^{\otimes,{\rm op}}. \]
Furthermore,
we have a functor
\[ {\rm Alg}_{\mathcal{P}\otimes\mathcal{O}}^{\rm Pair}
   (\catkappa) 
   \longrightarrow
   {\rm Mon}_{\mathcal{P}}({\rm Mon}_{\mathcal{O}}^{\rm oplax}(\cat)),\]
which associates to a pair $(\mathcal{C}^{\otimes},A)$
the mixed $(\mathcal{P},\mathcal{O})$-monoidal
$\infty$-category ${\rm Mod}_A^{\mathcal{O}}
(\mathcal{C})^{\otimes,\vee}
\to \mathcal{P}^{\boxtimes}
\times\mathcal{O}^{\otimes,{\rm op}}$.
\end{theorem}

We denote by ${\rm Pr}^L$
the large $\infty$-category of
presentable $\infty$-categories
and colimit-preserving functors.
By \cite[Proposition~4.8.1.15]{Lurie2},
the $\infty$-category ${\rm Pr}^L$ 
has a symmetric monoidal structure
in which the inclusion 
${\rm Pr}^L\hookrightarrow \wcat(\mathcal{K})$
is strong symmetric monoidal,
where $\mathcal{K}$ is the set of all small
simplicial sets.
By 
Theorem~\ref{thm:main-functor},
we obtain the following corollary.

\begin{corollary}\label{cor:presentable-duoidal-category}
Let $\mathcal{O}^{\otimes}$ be a small coherent $\infty$-operad 
and let $\mathcal{P}^{\boxtimes}$ be an $\infty$-operad. 
For a presentable $\mathcal{P}\otimes\mathcal{O}$-monoidal
$\infty$-category and a $\mathcal{P}\otimes\mathcal{O}$-algebra
object $A$ in $\mathcal{C}^{\otimes}$, 
we have a mixed $(\mathcal{P},\mathcal{O})$-monoidal
$\infty$-category
\[ {\rm Mod}_A^{\mathcal{O}}(\mathcal{C})^{\boxtimes, \otimes,\vee}
   \longrightarrow \mathcal{P}^{\boxtimes}\times
   \mathcal{O}^{\otimes,{\rm op}}. \]
Furthermore,
we have a functor
\[ {\rm Alg}_{\mathcal{P}\otimes\mathcal{O}}^{\rm Pair}
   ({\rm Pr}^L) 
   \longrightarrow
   {\rm Mon}_{\mathcal{P}}
   ({\rm Mon}_{\mathcal{O}}^{\rm oplax}(\wcat)),\]
which associates to a pair $(\mathcal{C}^{\otimes},A)$
the mixed $(\mathcal{P},\mathcal{O})$-monoidal
$\infty$-category 
${\rm Mod}_A^{\mathcal{O}}(\mathcal{C})^{\otimes,\vee}
   \to \mathcal{P}^{\boxtimes}\times
   \mathcal{O}^{\otimes,\rm op}$.
\end{corollary}

\section{$(\mathbb{E}_m,\mathbb{E}_n)$-duoidal
$\infty$-categories of $\mathbb{E}_n$-modules}
\label{section:em-en-duoidal}

Let $\mathbb{E}_k^{\otimes}$ be the little $k$-cubes operad 
for $k\ge 0$.
In this section we consider 
the important case in which  
$(\mathcal{P}^{\otimes},\mathcal{O}^{\otimes})
=(\mathbb{E}_m^{\otimes},\mathbb{E}_n^{\otimes})$
and $\mathcal{C}^{\otimes}$ is 
a presentable symmetric monoidal $\infty$-category.
The $\infty$-category ${\rm Mod}_A^{\mathbb{E}_n}(\mathcal{C})$
of $\mathbb{E}_n$-$A$-modules
has a structure of an $(\mathbb{E}_m,\mathbb{E}_n)$-duoidal
$\infty$-category
for any $A\in {\rm Alg}_{/\mathbb{E}_{m+n}}(\mathcal{C})$
by Corollary~\ref{cor:presentable-duoidal-category}.
The key to prove this
was lemma~\ref{lemma:f-shriek-boxproduct-compatibility}.
We will give
another proof of lemma~\ref{lemma:f-shriek-boxproduct-compatibility}
in this case
by using enveloping algebras.

We recall that $\mathbb{E}_k^{\otimes}$ is coherent
by \cite[Theorem~5.1.1.1]{Lurie2},
and that there is an
equivalence 
$\mathbb{E}_{m+n}^{\otimes}\simeq
\mathbb{E}_m^{\otimes}\otimes\mathbb{E}_n^{\otimes}$
for $m,n\ge 0$
by Dunn-Lurie Additivity Theorem
(\cite[Theorem~5.1.2.2]{Lurie2}).
In this section
we assume that $\mathcal{C}^{\otimes}$ is a presentable
symmetric monoidal $\infty$-category.
Under this assumption,
there is an equivalence 
${\rm Mod}_A^{\mathbb{E}_n}(\mathcal{C})\simeq
{\rm LMod}_{U(A)}(\mathcal{C})$
of $\infty$-categories
for $A\in {\rm Alg}_{\mathbb{E}_n}(\mathcal{C})$,
where ${\rm LMod}_{U(A)}(\mathcal{C})$
is the $\infty$-category of 
left modules over the enveloping algebra $U(A)$ of $A$.
Note that $U(A)$ is equivalent
to the factorization homology
$\int_{S^{n-1}}A$ 
by \cite[Proposition~3.16]{Francis}
(see also \cite[Example~5.5.4.16]{Lurie2}).

A map $f: A\to B$ in ${\rm Alg}_{\mathbb{E}_n}(\mathcal{C})$
induces a map
$U(A)\to U(B)$
of enveloping algebras
in ${\rm Alg}(\mathcal{C})$.
If we regard $M\in {\rm Mod}_A^{\mathbb{E}_n}(\mathcal{C})$
as a left $U(A)$-module,
then there is an equivalence
\[ f_!(M)\simeq U(B)\otimes_{U(A)}M \]
in ${\rm LMod}_{U(B)}(\mathcal{C})$.

Now, we recall that symmetric monoidal
structures on ${\rm Alg}_{\mathbb{E}_n}(\mathcal{C})$
and ${\rm Mod}^{\mathbb{E}_n}(\mathcal{C})$.
The $\infty$-category
${\rm Alg}_{\mathbb{E}_n}(\mathcal{C})$
is a symmetric monoidal $\infty$-category
by pointwise tensor product
(\cite[Example~3.2.4.4]{Lurie2}).
Furthermore,
${\rm Mod}^{\mathbb{E}_n}(\mathcal{C})$
is also a symmetric monoidal $\infty$-category
by pointwise tensor product 
such that the projection
${\rm Mod}^{\mathcal{O}}(\mathcal{C})
\to {\rm Alg}_{\mathcal{O}}(\mathcal{C})$
is a strong symmetric monoidal functor.

The following proposition
is a counterpart of Lemma~\ref{lemma:f-shriek-boxproduct-compatibility}
in the setting of this section.

\begin{proposition}
\label{prop:another-proof-En-envelop}
Let $\mathcal{C}^{\otimes}$ be a presentable symmetric monoidal
$\infty$-category.
Suppose that we have 
maps $f_i: A_i\to B_i$ in ${\rm Alg}_{\mathbb{E}_n}(\mathcal{C})$ 
and $M_i\in {\rm Mod}_{A_i}^{\mathbb{E}_n}(\mathcal{C})$
for $i=1,\ldots,n$.
Then the canonical map
\[ (f_1\otimes\cdots\otimes f_n)_!
   (M_1\otimes\cdots \otimes M_n)
   \longrightarrow
   (f_1)_!(M_1)\otimes\cdots\otimes (f_n)_!(M_n)\]
is an equivalence
for each $n\ge 0$.
\end{proposition}

\proof
It suffices to prove the case $n=2$.
We have an equivalence
\[ (f_1)_!(M_1)\otimes
   (f_2)_!(M_2)\simeq
   (U(B_1)\otimes U(B_2))\otimes_{(U(A_1)\otimes U(A_2))}
   (M_1\otimes M_2).\]
By \cite[Theorem~5.5.3.2]{Lurie2},
$U(A\otimes B)\simeq U(A)\otimes U(B)$
for any $A,B\in {\rm Alg}_{\mathbb{E}_n}(\mathcal{C})$. 
This implies an equivalence
\[ (U(B_1)\otimes U(B_2))\otimes_{(U(A_1)\otimes U(A_2))}
   (M_1\otimes M_2)\simeq
   (f_1\otimes f_2)_!(M_1\otimes M_2),\]
which completes the proof.
\qed

\bigskip

By using Proposition~\ref{prop:another-proof-En-envelop},
we obtain the following theorem
in the same way as in \S\ref{subsection:main-results}.

\begin{theorem}
\label{thm:Em-En-case}
Let $\mathcal{C}^{\otimes}$ be a presentable 
symmetric monoidal $\infty$-category.
For any $\mathbb{E}_{m+n}$-algebra
object $A$ in $\mathcal{C}^{\otimes}$,
we have a mixed $(\mathbb{E}_m,\mathbb{E}_n)$-monoidal
$\infty$-category
\[ {\rm Mod}_A^{\mathbb{E}_n}(\mathcal{C})^{\boxtimes,\otimes,\vee}
   \longrightarrow 
   \mathbb{E}_m^{\otimes}\times\mathbb{E}_n^{\otimes,\rm op},\]
in which the underlying $\infty$-category
is equivalent to
the $\infty$-category 
${\rm Mod}_A^{\mathbb{E}_n}(\mathcal{C})$
of $\mathbb{E}_n$-$A$-modules in $\mathcal{C}^{\otimes}$.
For any map $f: A\to B$ in ${\rm Alg}_{\mathbb{E}_n}(\mathcal{C})$,
the induced functor
$f_!: {\rm Mod}_A^{\mathcal{O}}(\mathcal{C})\to
{\rm Mod}_B^{\mathcal{O}}(\mathcal{C})$
is a bilax $(\mathbb{E}_m,\mathbb{E}_n)$-monoidal functor
that is strong monoidal
with respect to the $\mathbb{E}_m$-monoidal product.
\end{theorem}


\begin{thebibliography}{99}

\bibitem{Aguiar-Mahajan}
M. Aguiar and S. Mahajan, 
Monoidal functors, species and Hopf algebras,
CRM Monograph Series, 29. 
American Mathematical Society, Providence, RI, 2010. 






\bibitem{Francis}
J, Francis, 
The tangent complex and Hochschild cohomology of $\mathcal{E}_n$-rings. 
Compos. Math. 149 (2013), no. 3, 430--480. 


\bibitem{GHN}
D. Gepner, R. Haugseng, and T. Nikolaus, 
Lax colimits and free fibrations in $\infty$-categories,
preprint, 
arXiv:1501.02161v1.



\bibitem{Haugseng2}
R. Haugseng,
A fibrational mate correspondence for $\infty$-categories,
preprint,
arXiv:2011.08808.

\bibitem{HLN}
F. Hebestreit, S. Linskens, and J. Nuiten,
Orthofibrations and monoidal adjunctions, preprint,
arXiv:2011.11042.

\bibitem{Lurie1}
J. Lurie, 
Higher topos theory,
Annals of Mathematics Studies, 170. 
Princeton University Press, Princeton, NJ, 2009.

\bibitem{Lurie2}
J. Lurie,
Higher algebra,
available at 
http://www.math.harvard.edu/\~{}lurie/.



\bibitem{Torii1}
T. Torii,
On duoidal $\infty$-categories, 
preprint,
arXiv:2106.14121.

\bibitem{Torii3}
T. Torii,
On higher monoidal $\infty$-categories,
preprint,
arXiv:2111.00158.

\bibitem{Torii4}
T. Torii,
Perfect pairing for monoidal adjunctions,
preprint,
arXiv:2202.02493.

\end{thebibliography}
\end{document}